\documentclass[12pt]{amsart}
\usepackage{amsmath, amssymb, latexsym, amsthm, mathrsfs, color, mathtools, tikz,pgfplots,tikz-cd, graphicx, caption, stackengine,url,accents, scrextend, bm, mathrsfs}
\usepackage[top=2.5cm, bottom=2.5cm, left=2.5cm, right=2.5cm]{geometry}

\stackMath
\newcommand\tsup[2][2]{%
 \def\useanchorwidth{T}%
  \ifnum#1>1%
    \stackon[-1pt]{\tsup[\numexpr#1-1\relax]{#2}}{\hspace{1pt}\scriptstyle\sim}%
  \else%
    \stackon[.5pt]{#2}{\hspace{1pt}\scriptstyle\sim}%
  \fi%
}

\newcommand{\nc}{\newcommand}

\nc{\cO}{\mathcal{O}}

\newcommand{\sww}[1]{\mathsf{S}_1(\Omega,\Omega)}

\nc{\ram}{\mathsf{Ramsey}}
\nc{\soo}[1]{\mathsf{S}_1(\Op,\Op)}
\nc{\swl}[1]{\mathsf{S}_1(\Om,\Lambda)}
\nc{\swg}[1]{\mathsf{S}_1(\Om,\Ga)}
\nc{\goo}[1]{\gone(\Op,\Op)}
\nc{\gwo}{\gone(\Om,\cO)}
\nc{\gwl}[1]{\gone(\Om(#1),\Lambda(#1))}
\nc{\soox}[1]{\mathsf{S}_1(\cO(#1),\cO(#1))}
\nc{\swlx}[1]{\mathsf{S}_1(\Om(#1),\Lambda(#1))}
\nc{\goox}[1]{\gone(\cO(#1),\cO(#1))}
\nc{\gwox}[1]{\gone(\Om(#1),\cO(#1))}
\nc{\gwlx}[1]{\gone(\Om(#1),\Lambda(#1))}

\nc{\sfinwl}{\mathsf{S}_{\mathrm{fin}}(\Om,\Lambda)}
\nc{\sfinww}[1]{\mathsf{S}_{\mathrm{fin}}(\Om,\Om)}
\nc{\gfinoo}{\gfin(\cO,\cO)}
\nc{\gfinwo}{\gfin(\Om,\cO)}
\nc{\gfinwl}{\gfin(\Om,\Lambda)}
\nc{\sfinoox}[1]{\mathsf{S}_{\mathrm{fin}}(\cO(#1),\cO(#1))}
\nc{\sfinwlx}[1]{\mathsf{S}_{\mathrm{fin}}(\Om,\Lambda)}
\nc{\sfinwwx}[1]{\mathsf{S}_{\mathrm{fin}}(\Om(#1),\Om(#1))}
\nc{\gfinoox}[1]{\gfin(\cO(#1),\cO(#1))}
\nc{\gfinwox}[1]{\gfin(\Om(#1),\cO(#1))}
\nc{\gfinwlx}[1]{\gfin(\Om(#1),\Lambda(#1))}

\DeclareMathOperator{\C}{C_p}
\DeclareMathOperator{\fin}{Fin}

\nc{\mc}{\mathcal}
\nc{\thusfar}{\my{--- Edited thus far ---}}
\nc{\lei}{\le^\oo}
\nc{\sqsubs}{\sqsubseteq^*}
\nc{\card}[1]{\left|#1\right|}
\nc{\medcard}[1]{\biggl|\,#1\,\biggr|}
\nc{\smallcard}[1]{|\,#1\,|}
\nc{\bds}{bidirectional $\roth$-scale}
\nc{\bfP}{\mathbf{P}}
\nc{\bfQ}{\mathbf{Q}}
\nc{\bbT}{\mathbb{T}}
\nc{\bbZ}{\mathbb{Z}}
\nc{\bbN}{\mathbb{N}}
\nc{\bbC}{\mathbb{C}}
\nc{\beq}{\begin{equation}}\nc{\eeq}{\end{equation}}
\nc{\mbq}{\mb{?}}
\nc{\mb}[1]{{\mbox{\textbf{#1}}}}
\nc{\nop}{$\times$}
\nc{\fbn}{\!\!\fbox{\!\nop\!}\!\!}
\nc{\yup}{\checkmark}
\nc{\forces}{\Vdash}
\nc{\name}[1]{\dot{#1}}
\nc{\tf}{\my{FINISHED THUS FAR}}
\nc{\FU}{Fr\'echet--Urysohn}
\nc{\gs}{$\gamma$~space}
\nc{\Gab}{\Gamma_{\mathrm{B}}}
\nc{\Omb}{\Omega_{\mathrm{B}}}
\nc{\Ga}{\Gamma}
\nc{\Om}{\Omega}
\nc{\smallbinom}[2]{\begin{psmallmatrix} #1\\ #2 \end{psmallmatrix}}
\nc{\bgamma}{\smallbinom{\Om}{\Ga}}
\newcommand{\two}{\{0,1\}}
\nc{\productive}[2]{(#1,\allowbreak #2)^\x}
\nc{\prdct}[1]{#1^\x}
\nc{\Sel}{\mathsf{S}}
\nc{\sset}[2]{\{\,#1 : #2\,\}}
\nc{\smb}[1]{{\!\!\mb{#1}\!\!}}
\nc{\medset}[2]{{\biggl\{\,#1 : #2\,\biggr\}}}
\nc{\smallmedset}[2]{{\bigl\{\,#1 : #2\,\bigr\}}}
\nc{\set}[2]{{\left\{\,#1 : #2\,\right\}}}
\nc{\eseq}[1]{#1_1, \allowbreak #1_2, \allowbreak\dotsc} 

\nc{\eseqint}[3]{#1_{#2}, \allowbreak\dotsc,\allowbreak #1_{#3}} 
\nc{\eseqstart}[2]{#1_{#2},\allowbreak #1_{#2+1},\dotsc } 
\nc{\eprod}[1]{#1_{1}\times \allowbreak#1_{2}\times\dotsb}

\nc{\shortprod}[1]{\prod_{n=1}^\infty{#1}_n}

\nc{\eprodint}[3]{#1_{#2}\times \allowbreak\dotsb\times\allowbreak #1_{#3}}

\nc{\seleseq}[1]{#1_1\in \mathcal{#1}_1, \allowbreak #1_2\in \mathcal{#1}_2, \allowbreak\dotsc}
\nc{\cube}{(\Cantor)^\bbN}
\nc{\Match}{\op{Match}}
\nc{\concat}[1]{\hat{\phantom{a}}\langle #1\rangle}
\nc{\poset}{\mathbb{P}}
\nc{\fn}[1]{{\op{Fn}(#1\times\w,2)}}
\nc{\linadd}{\op{linadd}}
\nc{\nonprod}{\non^\x}
\nc{\alephes}{{\aleph_0}}
\nc{\my}[1]{{\color{red}{#1}}}
\nc{\later}[1]{{\color{green} #1}}
\nc{\BTs}[1]{{\color{green} #1 (BT)}}
\nc{\Cp}{\op{C}_\mathrm{p}}
\nc{\Bp}{\op{B}_p}
\nc{\Pa}[8]{\bibitem{#1} {#2}, \emph{#3}, {#4} \textbf{#5} ({#6}), {#7}--{#8}.}
\nc{\tPa}[5]{\bibitem{#1} {#2}, \emph{#3}, {#4}, to appear.}
\nc{\sPa}[4]{\bibitem{#1} {#2}, \emph{#3}, {#4}, submitted.}
\nc{\Bc}[9]{\bibitem{#1} {#2}, \emph{#3}, in: \textbf{#4} (#5), #6 #7, #8--#9.}
\nc{\fD}{\mathfrak{D}}
\nc{\fX}{\mathfrak{X}}
\nc{\Onbd}{\Op_{\mathrm{nbd}}} 
\nc{\Omnb}{\Om_{\mathrm{nbd}}} 
\nc{\od}{\mathfrak{od}}
\nc{\Setting}[7]{\xymatrix@R=4pt@C=7pt{#1\ar@{-}[r]&#2\ar@{-}[r]&#3\\&#4\ar@{-}[u]\\
#5\ar@{-}[uu]\ar@{-}[r] & #6\ar@{-}[u]\ar@{-}[r] & #7\ar@{-}[uu]}}
\nc{\mx}[1]{\begin{matrix}#1\end{matrix}}
\nc{\plim}{p\txt{-}\lim}
\nc{\Bgp}{{\Z^\bbN}}
\nc{\Cgp}{{{\Z_2}^\bbN}}
\nc{\Cite}[1]{\textbf{[#1]}}
\nc{\Next}[1]{{#1^+}}
\nc{\cFin}{\mathrm{cF}}
\nc{\scsp}{\text{-scale space}}
\nc{\cfn}{\text{cofinal}\ }
\nc{\Con}{\text{Concentrated}}
\nc{\Lind}{\text{Lindel\"of}\,}
\nc{\con}{\text{-Concentrated}}
\nc{\lind}{\text{-Lindel\"of}\,}
\nc{\ctbl}{\text{countably }\allowbreak}
\nc{\Hur}{\text{Hurewicz}}
\nc{\intvl}[2]{{[#1(#2),\allowbreak #1(#2\!+\!1))}}
\nc{\Bdd}{\mathbf{B}}
\nc{\Dfin}{\mathfrak{D}_\mathrm{fin}}
\nc{\grbl}{{\mbox{\textit{\tiny gp}}}}
\nc{\bbP}{\mathbb{P}}
\nc{\BOfat}{\B_{\Om_{\mathrm{fat}}}}
\nc{\Bgood}{\B_{\mathrm{good}}}
\nc{\compactN}{\cl{\mathbb{N}}}
\nc{\blocks}[2]{\op{cl}_{#2}(#1)}
\nc{\blocksplus}[2]{\op{cl}^+_{#2}(#1)}
\nc{\arx}[1]{\texttt{http://arxiv.org/math/#1}}
\nc{\bq}{\begin{quote}}
\nc{\eq}{\end{quote}}
\nc{\cl}[1]{\overline{#1}}
\nc{\Cl}[2]{\mathrm{cl}_{#1}(#2)}
\nc{\CH}{the Continuum Hypothesis}
\nc{\MA}{Martin's Axiom}
\nc{\Bfat}{\B_\mathrm{fat}}
\nc{\inv}{^{-1}}
\nc{\Cantor}{{\two^\bbN}}
\nc{\bP}{\mathbf{P}}
\nc{\bof}{\op{\fb}}
\nc{\dof}{\op{\fd}}
\nc{\bofF}{\bof(\cF)}
\nc{\sr}[3]{\underset{\mbox{#3}}{\mbox{#1}}}
\nc{\gp}{\binom{\Om}{\Ga}}
\nc{\gpsmall}{\mbox{$\gp$}}
\nc{\gig}{\gimel}
\nc{\gns}{\sone(\Om,\gig)}
\nc{\nsr}[2]{#1}
\nc{\Srg}{{\mathbb{S}}}
\nc{\Srgs}{{\mathbb{S}^*}}
\nc{\NN}{{\bbN^{\bbN}}}
\nc{\ZN}{{\Z^{\bbN}}}
\nc{\NNup}{{\bbN^{\uparrow\bbN}}}
\nc{\NNupb}{{b^{\uparrow\bbN}}}
\nc{\Pof}{\op{P}}
\nc{\PN}{{\Pof(\bbN)}}
\nc{\rothx}[1]{{[#1]^{\mbox{\tiny $\infty$}}}}
\nc{\tx}{{\tilde{x}}}
\nc{\roth}{{[\bbN]^{\mbox{\tiny $\infty$}}}} 
\nc{\roths}{{[b]^{\mbox{\tiny $\infty$}}}} 
\nc{\Fin}{\mathrm{Fin}(\bbN)}
\nc{\ici}{[\bbN]^{ \infty, \infty}}
\nc{\Inc}{{\compactN^{\uparrow\bbN}}}
\nc{\powInc}[1]{{\big(\Inc\big)^{#1}}}
\nc{\powFin}[1]{{\big(\Fin\big)^{#1}}}
\nc{\powPN}[1]{{\big(\PN\big)^{#1}}}
\nc{\NcompactN}{{\compactN^\bbN}}
\nc{\seq}[1]{\la #1\ra_{n\in\bbN}}
\nc{\Uarrow}{\smash{\big\uparrow}}
\nc{\LE}{\preccurlyeq}
\nc{\GE}{\succcurlyeq}
\nc{\op}{\operatorname}
\nc{\im}{\op{Im}}
\nc{\Span}{\op{span}}
\nc{\maxfin}{\op{maxfin}}
\nc{\ran}{\op{range}}
\nc{\iso}{\cong}
\nc{\Madd}{{\M}^\star}
\nc{\cI}{\mathcal{I}}
\nc{\cJ}{\mathcal{J}}
\nc{\scrA}{\mathscr{A}}
\nc{\scrB}{\mathscr{B}}
\nc{\scrC}{\mathscr{C}}
\nc{\scrD}{\mathscr{D}}
\nc{\scrF}{\mathscr{F}}
\nc{\scrK}{\mathscr{K}}
\nc{\A}{\D\forall}
\nc{\B}{\mathrm{B}}
\nc{\cB}{\mathcal{B}}
\nc{\cZ}{\mathcal{Z}}
\nc{\bB}{\mathbf{B}}
\nc{\BS}{\mathbf{B}(\mathcal{S})}
\nc{\BF}{\mathbf{B}(\mathcal{F})}
\nc{\BU}{\mathbf{B}(\mathcal{U})}
\nc{\cSp}{\mathcal{S}^+}
\nc{\cFp}{\mathcal{F}^+}
\nc{\cUp}{\mathcal{U}^+}

\nc{\BG}{\B_\Ga}
\nc{\BL}{\B_\Lambda}
\nc{\BT}{\B_\Tau}
\nc{\BTstar}{\B_{\Tau^*}}
\nc{\BO}{\B_\Om}
\nc{\DO}{\cD_\Om}
\nc{\KO}{\cK_\Om}
\nc{\CG}{C_\Ga}
\nc{\CL}{C_\Lambda}
\nc{\CT}{C_\Tau}
\nc{\CTstar}{C_{\Tau^*}}
\nc{\CO}{C_\Om}
\nc{\COgp}{C_{\Om^{\grbl}}}
\nc{\CLgp}{C_{\Lambda^{\grbl}}}
\nc{\BOgp}{\B_{\Om}^{\grbl}}
\nc{\BLgp}{\B_{\Lambda^{\grbl}}}
\nc{\sfC}{\mathsf{C}}
\nc{\sfD}{\mathsf{D}}
\nc{\bD}{\mathbf{D}}
\nc{\Tau}{\mathrm{T}}
\nc{\cA}{\mathcal{A}}
\nc{\cK}{\mathcal{K}}
\nc{\cD}{\mathcal{D}}
\nc{\cE}{\mathcal{E}}
\nc{\cF}{\mathcal{F}}
\nc{\cS}{\mathcal{S}}
\nc{\cT}{\mathcal{T}}
\nc{\cG}{\mathcal{G}}
\nc{\cY}{\mathcal{Y}}
\nc{\J}{\mathcal{J}}
\nc{\cL}{\mathcal{L}}
\nc{\cM}{\mathcal{M}}
\nc{\cN}{\mathcal{N}}
\nc{\cH}{\mathcal{H}}

\nc{\scF}{\mathscr{F}}
\nc{\scH}{\mathscr{H}}

\nc{\Op}{\mathrm{O}}
\nc{\rmA}{\mathrm{A}}
\nc{\rmF}{\mathrm{F}}
\nc{\rmB}{\mathrm{B}}
\nc{\rmD}{\mathrm{D}}
\nc{\rmP}{\mathrm{P}}
\nc{\cC}{\mathcal{C}}
\nc{\cP}{\mathcal{P}}
\nc{\bbQ}{\mathbb{Q}}
\nc{\bbR}{\mathbb{R}}
\nc{\cU}{\mathcal{U}}
\nc{\cQ}{\mathcal{Q}}
\nc{\Un}{\bigcup}
\nc{\cV}{\mathcal{V}}
\nc{\cR}{\mathcal{R}}
\nc{\tcR}{\tilde{\mathcal{R}}}
\nc{\cW}{\mathcal{W}}
\nc{\Z}{{\mathbb Z}}
\nc{\Impl}{\Rightarrow}
\long\def\forget#1\forgotten{\marginpar{\textcolor{green}{Forgetting...}}}
\nc{\ft}{\mathfrak{t}}
\nc{\fb}{\mathfrak{b}}
\nc{\fc}{\mathfrak{c}}
\nc{\fd}{\mathfrak{d}}
\nc{\fg}{\mathfrak{g}}
\nc{\oo}{\infty}
\nc{\fr}{\mathfrak{r}}
\nc{\fk}{\mathfrak{k}}
\nc{\bidi}{\mathfrak{bidi}}
\nc{\fu}{\mathfrak{u}}
\nc{\fh}{\mathfrak{h}}
\nc{\fp}{\mathfrak{p}}
\nc{\fj}{\mathfrak{j}}
\nc{\fs}{\mathfrak{s}}
\nc{\w}{\omega}
\nc{\x}{\times}
\nc{\Iff}{\Leftrightarrow}

\nc{\nin}{\notin}
\nc{\cat}{\hat{\ }}
\nc{\sub}{\subseteq}
\nc{\spst}{\supseteq}
\nc{\sm}{\setminus}
\nc{\as}{\subseteq^*}
\nc{\les}{\le^*}
\nc{\leinf}{\le^{\infty}}
\nc{\leS}{\le_S}
\nc{\leF}{\le_F}
\nc{\leU}{\le_U}
\nc{\rest}{\restriction}

\nc{\la}{\langle}
\nc{\ra}{\rangle}
\nc{\E}{\exists}
\nc{\dom}{\op{dom}}
\nc{\cov}{\op{cov}}
\nc{\add}{\op{add}}
\nc{\addmen}{\add(\Men{})}
\nc{\cof}{\op{cof}}
\nc{\cf}{\op{cf}}
\nc{\non}{\op{non}}
\nc{\unif}{\op{non}}
\nc{\COV}{\op{COV}}
\nc{\ADD}{\op{ADD}}
\nc{\COF}{\op{COF}}
\nc{\NON}{\op{NON}}
\nc{\impl}{\to}
\nc{\Lp}{\mathcal{L_\p}}
\nc{\Wlog}{without loss of generality}
\newtheorem{thm}{Theorem}[section]
\nc{\bthm}{\begin{thm}} \nc{\ethm}{\end{thm}}
\newtheorem{need}[thm]{Need}
\nc{\bneed}{\begin{need}\color{dg}} \nc{\eneed}{\end{need}}
\newtheorem{prop}[thm]{Proposition}
\nc{\bprp}{\begin{prop}} \nc{\eprp}{\end{prop}}
\newtheorem{fact}[thm]{Fact}
\nc{\bfct}{\begin{fact}} \nc{\efct}{\end{fact}}
\newtheorem{prob}[thm]{Problem}
\nc{\bprb}{\begin{prob}} \nc{\eprb}{\end{prob}}
\newtheorem{lem}[thm]{Lemma}
\nc{\blem}{\begin{lem}} \nc{\elem}{\end{lem}}
\newtheorem{app}[thm]{Application}
\nc{\bapp}{\begin{app}} \nc{\eapp}{\end{app}}
\newtheorem{claim}[thm]{Claim}
\nc{\bclm}{\begin{claim}} \nc{\eclm}{\end{claim}}
\newtheorem{cor}[thm]{Corollary}
\nc{\bcor}{\begin{cor}} \nc{\ecor}{\end{cor}}
\newtheorem{conj}[thm]{Conjecture}
\nc{\bcnj}{\begin{conj}} \nc{\ecnj}{\end{conj}}
\theoremstyle{definition}
\newtheorem{defn}[thm]{Definition}
\nc{\bdfn}{\begin{defn}} \nc{\edfn}{\end{defn}}
\newtheorem{obs}[thm]{Observation}
\nc{\bobs}{\begin{obs}} \nc{\eobs}{\end{obs}}
\theoremstyle{remark}
\newtheorem{rem}[thm]{Remark}
\nc{\brem}{\begin{rem}} \nc{\erem}{\end{rem}}
\newtheorem{cnv}[thm]{Convention}
\nc{\bcnv}{\begin{cnv}} \nc{\ecnv}{\end{cnv}}
\newtheorem{exam}[thm]{Example}
\nc{\bexm}{\begin{exam}} \nc{\eexm}{\end{exam}}
\nc{\bpf}{\begin{proof}} \nc{\epf}{\end{proof}
}
\nc{\be}{\begin{enumerate}}
\nc{\ee}{\end{enumerate}}
\nc{\bi}{\begin{itemize}}
\nc{\bimy}{\my{\begin{itemize}}
\nc{\eimy}{\end{itemize}}}
\nc{\itm}{\item}
\nc{\ei}{\end{itemize}}
\nc{\Subsection}[1]{\goodbreak\subsection*{#1}}

\nc{\sone}{\mathsf{S}_1}
\nc{\sfin}{\mathsf{S}_\mathrm{fin}}
\nc{\ufin}{\mathsf{U}_\mathrm{fin}}
\nc{\Split}{\mathsf{Split}}

\nc{\gone}{\mathsf{G}_1}    
\nc{\tgfin}{{\mathsf{G}}^*_\mathrm{fin}}
\nc{\gfin}{\mathsf{G}_\mathrm{fin}}

\nc{\men}[1]{\sfin(\Op(#1),\Op(#1))}
\nc{\sch}{\ufin(\cO,\Omega)}
\nc{\rothb}{\text{Rothberger}}
\nc{\pmen}{\sfin(\Omega,\Omega)}
\nc{\Rothb}{\sone(\Op,\Op)}

\nc{\prothb}{\sone(\Omega,\Omega)}
\nc{\tU}{{\tilde{U}}}
\nc{\tF}{{\tilde{F}}}
\nc{\tY}{{\tilde{Y}}}
\nc{\tX}{{\tsup[1]{X}}}
\nc{\dtX}{{\tsup[2]{X}}}
\nc{\dt}[1]{{\tsup[2]{#1}}}
\nc{\td}{{\tilde{d}}}
\nc{\tz}{{\tilde{z}}}
\nc{\cfd}{\cf(\fd)}

\nc{\msep}{\sfin(\cD,\cD)}
\nc{\rsep}{\sone(\cD,\cD)}
\nc{\cft}{\sfin(\Omega_{\mathbf{0}},\Omega_{\mathbf{0}})}
\nc{\scft}{\sone(\Omega_{\mathbf{0}},\Omega_{\mathbf{0}})}

\nc{\Umen}{U\text{-Menger}}
\nc{\hur}{\ufin(\cO,\Gamma)}
\nc{\tUmen}{\tU\text{-Menger}}

\nc{\Men}{\text{Menger}}
\nc{\Sch}{\text{Scheepers}}

\nc{\aspst}{\prescript{*}{}{\spst}\ }
\nc{\eqs}{=^*}
\nc{\ctblOm}{\Omega_{\mathrm{ctbl}}}
\nc{\GNga}{{\smallbinom{\Om}{\Ga}}}
\nc{\ctblga}{\smallbinom{\ctblOm}{\Ga}}

\nc{\nadd}{\cN_{\mathrm{add}}}
\nc{\ball}{\mathrm{B}}
\nc{\cOunif}{\cO^{\textrm{unif}}}

\nc{\sep}{
\vspace{2cm}
\noindent
\begin{minipage}{\textwidth}
	\textcolor{red}{\rule{\textwidth}{1pt}}
\end{minipage}
}
\nc{\FS}{\op{FS}}
\nc{\sums}{\op{SS}}
\nc{\SG}{\op{SG}}
\nc{\tSG}{\SG_\odot}
\nc{\G}{\op{G}}
\nc{\FSG}{\op{FSG}}
\nc{\FP}{\op{FP}}
\nc{\nonNadd}{\non(\nadd)}
\nc{\borga}{\Ga_\mathrm{Bor}}
\nc{\pick}{x}
\nc{\gen}{y}
\nc{\nullzind}{\sone(\{\Op_n^{\mathsf{unif}}\}_{n\in\bbN},\Ga)}
\nc{\nullzindf}[1]{\sone(\{\Op_{#1}^{\mathsf{unif}}\}_{n\in\bbN},\Ga)}

\definecolor{dg}{RGB}{42,101,24}

\nc{\myb}[1]{\textcolor{blue}{#1}}
\nc{\mydg}[1]{{\color{dg}{#1}}}
\DeclareMathOperator{\eexists}{\exists}
\DeclareMathOperator{\fforall}{\forall}
\nc{\Exists}[1]{\bigl(\eexists #1\bigr)}
\nc{\Forall}[1]{\bigl(\fforall #1\bigr)}
\nc{\End}[1]{\bigl(#1\bigr)}
\nc{\dmo}[2]{\DeclareMathOperator{#1}{#2}}
\dmo{\Asc}{Asc}
\nc{\plusmin}{\wedge}

\nc{\cBsub}{{\cB^{\mbox{\tiny $\sub$}}}}

\nc{\Alice}{{\textsc{Alice}{}}}
\nc{\Bob}{{\textsc{Bob}}}

\nc{\bfO}{\mathbf{0}}

\nc{\proba}[1]{F}

\nc{\opn}[1]{\Op}
\nc{\sfinoo}[1]{\sfin(\Op,\Op)}
\nc{\om}[1]{\Om}

\usepackage[all]{xy}

\newlength{\spacebox}
\usetikzlibrary{arrows}

\title{}

\author[P.~Szewczak]{Piotr Szewczak}
\address{Piotr Szewczak, Institute of Mathematics, Faculty of Mathematics and Natural Science College of Sciences, Cardinal Stefan Wyszy\'nski University in Warsaw, W\'oycickiego 1$\slash$3, 01--938 Warsaw, Poland
}
\email{p.szewczak@wp.pl}
\urladdr{http://piotrszewczak.pl}

\makeatletter
\@namedef{subjclassname@2020}{\textup{2020} Mathematics Subject Classification}
\makeatother

\subjclass[2020]{Primary: 05D10, 54D20, 54C35; Secondary: 16W22}
		
\keywords{
finite colorings, semigroups, \v{C}ech--Stone compactification, combinatorial covering proerties, local properties in function spaces, infinite topological games, Menger's property, Rothberger's property, selection principles.
}

\title{Abstract colorings, games and ultrafilters}

\begin{document}

\maketitle

\centerline{\emph{Dedicated to Boaz Tsaban}}

\begin{abstract}

The main result provide a common generalization for Ramsey-type theorems concerning finite colorings of edge sets of complete graphs with vertices in infinite semigroups.
We capture the essence of theorems proved in different fields: for natural numbers due to Milliken--Tylor, Deuber--Hindman, Bergelson--Hindman, for combinatorial covering properties due to Scheepers and Tsaban, and local properties in function spaces due to Scheepers.
To this end, we use idempotent ultrafilters in the \v{C}ech--Stone compactifications of discrete infinite semigroups and topological games.
The research is motivated by the recent breakthrough work of Tsaban about colorings and the Menger covering property.

\end{abstract}

\section{Background}

\subsection{Colorings and natural numbers}

A \emph{coloring} of a nonempty set $X$ is a function $\chi\colon X\to \{1,\dotsc,k\}$, where $k$ is a natural number. 
Given a coloring $\chi$ of a set $X$, a set $A\sub X$ is  \emph{$\chi$-monochromatic} (or just \emph{monochromatic}, when $\chi$ is clear from the context), if there is a color $i$ such that $\chi(a)=i$ for all $a\in A$.
By the van der Waerden Theorem~\cite{vdW}, for each coloring of the set of natural numbers $\bbN$, there is a monochromatic arithmetic progression of an arbitrarily finite length.
In the comprehensive work, Bergelson and Hindman~\cite{BHcell} considered families of finite subsets of $\bbN$ with the property that for each coloring of $\bbN$, there is  a monochromatic set in the family.
By the result of Hindman~\cite[Theorem~6.7]{Hu}, a family of finite subsets of $\bbN$ has the above property if and only if there is an ultrafilter on $\bbN$ such that each set in the ultrafilter contains a set in the family.
Ultrafilters play an important role in consideration of colorings of $\bbN$, especially, when an algebraic structure of $\bbN$ is involved.

Let $S$ be an infinite semigroup with the discrete topology. Usually, we denote the semigroup operation by $+$, even if this operation is not commutative.
The \emph{\v{C}ech--Stone compactification}, $\beta S$, of $S$ is the family of all ultrafilters on $S$ with the topology generated by the sets $\sset{p\in\beta S}{A\in p}$, where $A\sub S$.
The operation $+$ on $S$ can be extended to the operation $+$ on $\beta S$ such that for ultrafilters $p,q\in\beta S$, we have
\[
A\in p+q\quad\text{if and only if}\quad \set{b\in S}{\Exists{C\in q}\End{b+C\sub A}}\in p.
\]
Then, $\beta S$ is a compact semigroup.
For more details about properties of $\beta S$, we refer to the book of Hindman and Strauss~\cite{HS}.
Central in our investigations is the Ellis--Numakura Lemma~\cite{E,N} which asserts that each compact nonempty subsemigroup of $\beta S$ has an \emph{idempotent}, i.e., an ultrafilter $e\in\beta S$ with $e+e=e$.

By the celebrated Hindman Finite Sums Theorem~\cite{HFS}, for each coloring of $\bbN$ there is an infinite subset of $\bbN$ such that all finite sums of its elements (in usual sense) have the same color.
Galvin and Glazer provided an elegant and short proof of the Hindman Theorem, based on the Ellis--Numakura Lemma.

Let $\cA,\cR$ be families of nonempty sets.
The family $\cA$ is \emph{large} for the family $\cR$ if any set in $\cA$ contains a set in $\cR$.
The family $\cA$ is \emph{large for a sequence $\eseq{\cR}$} of families of nonempty sets if it is large for all families $\cR_n$.
For a nonempty set $A$, let $\rothx{A}$ be the set of all \emph{infinite} subsets of $A$ and $\fin(A)$ be  the family of all \emph{finite nonempty} subsets of $A$.

The Hindman Finite Sums Theorem is a special case of the following later theorem due to Deuber and Hindman. Originally, Theorem~\ref{thm:DH} was formulated in a slightly different way for specific families $\cR_n$; see Subsection~\ref{ssec:nat} for more details. 
In the below results of this Subsection, we consider $\bbN$ as a semigroup with the standard addition, as a semigroup operation. 

\bthm[{Deuber, Hindman~\cite{DH}}]
\label{thm:DH}
Assume that $\roth$ contains an idempotent, large for a sequence $\eseq{\cR}\sub\Fin$.
Then for each coloring of  $\bbN$, there are pairwise disjoint sets $R_1\in\cR_1, R_2\in\cR_2,\dotsc$ such that all finite sums choosing at most one from each set from the sequence $\eseq{R}$,
have the same color.
\ethm

In order to consider colorings in higher dimensions we need the following notations.
For a set $A$, let $[A]^2$ be the family of all two-element subsets of $A$, equivalently, the edge set of the complete graph with vertices in $A$.
For sets $H_1,H_2\in\Fin$, we write $H_1<H_2$ if $\max H_1<\min H_2$.

\bdfn
Let $\eseq{a}$ be a sequence in a semigroup $S$.
For a set $H=\{i_1,\dotsc,i_n\}\in\Fin$ with $i_1<\dotsb<i_n$, where $n$ is a natural number, let
\[
a_H:=a_{i_1}+\dotsb+a_{i_n}.
\]
A sequence $\eseq{a}$ is \emph{proper} if $a_{H_1}\neq a_{H_2}$ for all sets $H_1,H_2\in\Fin$ with $H_1<H_2$.
A \emph{sumgraph} of a proper sequence $\eseq{a}$ is the set 
\[
\smallmedset{\{a_{H_1}, a_{H_2}\}}{H_1,H_2\in\Fin\text{ and }H_1<H_2}.
\]
\edfn

\bthm[{Milliken--Taylor~\cite{mill,tylor}}]
For each coloring of  $[\bbN]^2$, there is a proper sequence $\eseq{a}$ in $\bbN$ with a monochromatic sumgraph.
\ethm

Let $\eseq{a}$ be a sequence in a semigroup $S$ and $n\leq m$ be natural numbers.
Let 
\[
\FS(\eseqint{a}{n}{m}):=\sset{a_H}{H\in\fin(\{n,\dotsc,m\})}.
\]
Bergelson and Hindman proved the following result, a common generalization for the Ramsey Theorem and the van der Waerden Theorem.

\bthm[{Bergelson, Hindman~\cite[Theorem~2.5.]{BH}}]
\label{thm:BHweak}
Assume that $m$ is a natural number and $\roth$ contains an idempotent, large for a sequence $\eseq{\cR}\sub\Fin$.
Then for each coloring of  $[\bbN]^m$, there are increasing sequences $x^{(n)}_1,x^{(n)}_2,\dotsc$ of natural numbers and pairwise disjoint sets $R_1\in\cR_1, R_2\in\cR_2,\dotsc$ with the following property.
For every natural numbers $n_1<n_2$ and elements
\begin{align*}
a_1&\in R_{n_1}\cup\Un\sset{\FS(x^{(l)}_{1},\dotsc,x^{(l)}_{n_1})}{l\leq n_1},\\
a_2&\in R_{n_2}\cup\Un\sset{\FS(x^{(l)}_{n_1+1},\dotsc,x^{(l)}_{n_2})}{n_1<l\leq n_2}
\end{align*}
with $a_1<a_2$, the sets $\{a_1,a_2\}$ have the same color.
\ethm

The above theorems, also in its higher-dimensional versions, are consequences of the main results of the paper, see Subsection~\ref{ssec:nat}.
Subsection~\ref{ssec:BH} is devoted to a modification of Theorem~\ref{thm:BHweak}, where both operations, additions and multiplications on $\bbN$, are involved.

\subsection{Covers}
Let $\cA,\cB$ be families of sets.
Define
\smallskip

\begin{labeling}{$\sfin(\cA,\cB)$:}
	\item [$\sfin(\cA,\cB)$:] for each sequence $\eseq{A}\in\cA$, there are finite sets $F_1\sub A_1,F_2\sub A_2,\ldots$ such that $\Un_{n\in\bbN} F_n\in \cB$.	
\end{labeling}
\smallskip

By \emph{space} we mean an infinite Tychonoff topological space.
A \emph{cover} of a space is a family of nonempty proper subsets of the space whose union is the entire space. 
Let $\opn{X}$ be the family of all \emph{open covers} of a space.
One of the central properties of spaces, in the \emph{selection principles theory}, is the \emph{Menger} property $\sfinoo{X}$~\cite{Menger24, Hure25}.
This property generalizes $\sigma$-compactness and it implies Lindel\"ofness.
In ZFC, there is a Menger \emph{set of reals}, i.e., a space homeomorphic with a subspace of the real line, which is not $\sigma$-compact~\cite{FrMill, BaTs}.
The Menger property has significant applications in other fields.
It characterizes filters whose Mathias forcing notion does not add a dominating function~\cite{ChRZd}.
One of the major open problems within set-theoretic topology is the \emph{D-space problem}, whether every Hausdorff Lindel\"of space is a D-space.
Thus far, the class of Menger spaces is the wider natural class of spaces for which a positive answer to the D-space problem is known~\cite{Aurichi}. 

The Menger property has also connections with local properties in function spaces.
A space $Y$ has \emph{countable fan tightness} if for every point $y\in Y$ and a sequence $\eseq{A}\sub Y$ such that $y\in\bigcap_{n\in\bbN}\overline{A_n}$, there are finite sets $F_1\sub A_1, F_2\sub A_2,\dotsc$ such that $y\in\overline{\Un_{n\in\bbN}F_n}$.
A cover of a space is an \emph{$\w$-cover} if each finite subset of the space is contained in some set from the cover.
Let $\om{X}$ be the family of all \emph{open $\w$-covers} of a space.
Let $X$ be a space and $\Cp(X)$ be the set of all continuous real-valued functions on $X$ with the pointwise convergence topology. 
The space $\Cp(X)$ has countable fan tightness if and only if $X$ is Menger in all finite powers~\cite{arhcft}, equivalently, the space $X$ satisfies $\sfinww{X}$~\cite[Theorem~3.9]{coc2}.
The statement that for sets of reals, the properties $\sfinoo{X}$ and $\sfinww{X}$ are equivalent, is independent from ZFC~\cite{Mill, Comb, pMWien}.

Tsaban observed that the Hindman Finite Sums Theorem and the Milliken--Tylor Theorem can be viewed as coloring theorems concerning countable covers of countable, discrete space.
Every countable space is Menger, and thus these theorems are special instances~\cite[Example~4.7]{tsramsey} of the following result.
For a space $X$, let $\tau$ be the topology of $X$.

\bthm[{Tsaban~\cite[Theorem~4.6]{tsramsey}}]
\label{thm:ts}
Let $X$ be a space satisfying the Menger property $\sfinoo{X}$ and $\cup$ be the semigroup operation on $\tau$.
Then for every infinite family $\cU\in \opn{X}$ and coloring of  $[\tau]^2$ there are pairwise disjoint sets $\eseq{\cF}\in\fin(\cU)$ with the following properties.
\be
\item The family $\Un_{n\in\bbN}\cF_n$ is in $\opn{X}$.
\item The sequence $\eseq{\Un\cF}$ is proper.
\item The sumgraph of the sequence $\eseq{\Un\cF}$ is monochromatic.
\ee
\ethm

Scheepers used colorings, in abstract context, to characterize combinatorial covering properties~\cite{coc1,schpart,coc6}.
Below, we present one of these applications.

A cover of a space is a \emph{large cover}, if each element of the space belongs to infinitely many sets in the cover.
For a space, let $\Lambda$ be the family of all open large covers of the space.

\bdfn
Let $S$ be an infinite semigroup.
A \emph{partite graph} of a sequence $\eseq{F}\in\fin(S)$ of pairwise disjoint sets, is the set
\[
\smallmedset{\{a_{i},a_{j}\}}{a_{i}\in F_{i}, a_{j}\in F_{j}\text{ and }i<j}.
\]
\edfn

\bthm[{Scheepers~\cite[Theorem~6]{schpart}}]
\label{thm:schsfinww}
A separable metrizable space $X$ satisfies the Menger property $\sfinoo{X}$ if and only if for every family $\cU\in\om{X}$ and a coloring of $[\tau]^2$, there are pairwise disjoint sets $\eseq{\cF}\in\fin(\cU)$ such that the family $\Un_{n\in\bbN}\cF_n$ is in $\Lambda$ and the partite graph of the sequence $\eseq{\cF}$ is monochromatic. 
\ethm

Other results of Scheepers, also about local properties of spaces, are discussed in details in Subsections~\ref{ssec:covers} and~\ref{ssec:loc}.
In Subsection~\ref{ssec:ts}, we combine Theorems~\ref{thm:ts} and~\ref{thm:schsfinww}.

\subsection{Games}
The above results of Tsaban and Scheepers, involve game-theoretical characterizations of covering properties.
Also in our considerations, we use this tool extensively.

\bdfn
Let $\cA$ and $\cB$ be nonempty families of sets.
A game $\gfin(\cA,\cB)$ is a game with two players \Alice{} and \Bob{}.
In the first inning, \Alice{} chooses a set $A_1\in \cA$ and \Bob{} replies with a finite set $F_1\sub A_1$.
In the second inning, \Alice{} chooses a set $A_2\in\cA$ and \Bob{} replies with a finite set $F_2\sub A_2$, etc.
\Bob{} wins, if the set $\Un_{n\in\bbN}F_n$ is in $\cB$, otherwise, \Alice{} wins.
\edfn

Let $\cA$, $\cB$ be nonempty families of sets. 
If \Alice{} has no winning strategy in the game $\gfin(\cA,\cB)$, then the statement $\sfin(\cA,\cB)$ holds. 
In many cases, this implication can be reversed, and thus some topological properties can be characterized in the language of games and winning strategies.
Among others, the most celebrated result is the Hurewicz Theorem.

\bthm[{Hurewicz~\cite[Theorem~10]{Hure25}}]
\label{thm:hur}
A space $X$ satisfies the Menger property $\sfinoo{X}$ if and only if \Alice{} has no winning strategy in the game $\gfin(\opn{X},\opn{X})$.
\ethm

A conceptual proof of Theorem~\ref{thm:hur} can be found in one of the previous work~\cite{Mgame}.
For a separable metrizable space, 
Bob has a winning strategy in the game $\gfin(\opn{X},\opn{X})$ if and only if the space is $\sigma$-compact (\cite[Corollary~4]{telg},~\cite{schtelg}).
For sets of reals, the game $\gfin(\opn{X},\opn{X})$ is undetermined, i.e., there is a space such that none of the players has a winning strategy in the game $\sfinoo{X}$.

\section{The main result}

The results concerning colorings from the introduction have higher-dimensional versions.
Thus, we need to extent above definitions to higher dimensions.

\bdfn
Let $S$ be an infinite semigroup and $m$ be a natural number.

\be 
\item An \emph{$m$-sumgraph} of a proper sequence $\eseq{a}\in S$ is the set 
\[
\SG^m(\eseq{a}):=
\smallmedset{\{a_{H_1},\dotsc, a_{H_m}\}}{H_1,\dotsc,H_m\in\Fin\text{ and }H_1<\dotsb<H_m}.
\]

\item A \emph{partite $m$-sumgraph} of a sequence $\eseq{F}\in\fin(S)$, where all sequences in $\eprod{F}$ are proper is the set
\[
\Un\smallmedset{\SG^m(\eseq{a})}{(\eseq{a})\in\eprod{F}}.
\]

\item A \emph{partite $m$-graph} of a sequence $\eseq{F}\in\fin(S)$ of pairwise disjoint sets, is the set
\[
\smallmedset{\{a_{i_1},\dotsc, a_{i_m}\}}{a_{i_1}\in F_{i_1},\dotsc,a_{i_m}\in F_{i_m}\text{ and }i_1<\dotsb<i_m}.
\]
\ee
\edfn

For a set $A$ and a natural number $m$, let $[A]^m$ be the family of all $m$-element subsets of $A$. For nonempty families $\cA,\cB$ of nonempty sets, a game $\gone(\cA,\cB)$ is the following game with two players \Alice{} and \Bob{}.
In the first inning, \Alice{} chooses a set $A_1\in \cA$ and \Bob{} replies with an element $a_1\in A_1$.
In the second inning, \Alice{} chooses a set $A_2\in\cA$ and \Bob{} replies with an element $a_2\in A_2$, etc.
\Bob{} wins, if $\sset{a_n}{n\in\bbN}\in\cB$, otherwise, \Alice{} wins.

\bthm\label{thm:mainfin}
Let $S$ be an infinite semigroup and $m$ be a natural number.
Assume that a family $\cA\sub\rothx{S}$ contains an idempotent, large for a sequence $\eseq{\cR}\sub\fin(S)$, and $\cB\sub\rothx{S}$ is a family such that \Alice{} has no winning strategy in the game $\gfin(\cA,\cB)$.
Then for each coloring of $[S]^m$, there are finite families $\cF_1\sub\cR_1, \cF_2\sub\cR_2,\dotsc$ of pairwise disjoint sets and finite sets $F_1\sub\Un\cF_1, F_2\sub\Un\cF_2,\dotsc$ with the following properties.
\be
\item The set $\Un_{n\in\bbN}F_n$ is in $\cB$.
\item All sequences in the product $\eprod{\Un\cF}$ are proper.
\item The partite $m$-sumgraph of the sequence $\eseq{\Un\cF}$ is monochromatic.
\ee 

Moreover, if \Alice{} has no winning strategy in the game $\gone(\cA,\cB)$, then we may require that the above families $\eseq{\cF}$ and sets $\eseq{F}$ are singletons, i.e., there are sets $R_1\in\cR_1, R_2\in\cR_2,\dotsc$ and elements $a_1\in R_1, a_2\in R_2,\dotsc$ with the following properties.
\be
\renewcommand{\theenumi}{\arabic{enumi}'}
\item The set $\sset{a_n}{n\in\bbN}$ is in $\cB$.
\item All sequences in the product $\eprod{R}$ are proper.
\item The $m$-sumgraph of the sequence $\eseq{R}$ is monochromatic.
\ee 
\ethm

In order to prove Theorem~\ref{thm:mainfin}, we need the following notions and auxiliary results.

\blem\label{lem:ult}
Let $S$ be an infinite set and $p\sub\rothx{S}$ be an ultrafilter, large for a family $\cR\sub\fin(S)$.
For each set $D\in p$, there is a family $\cR'\sub\cR$ of pairwise disjoint sets such that $\Un\cR'\sub D$ and $\Un\cR'\in p$.
\elem

\bpf
Let $\cP$ be the collection of all subfamilies of $\cR$ of pairwise disjoint sets, whose unions are contained in $D$, with a partial order $\sub$.
Since $A\in p$ and $p$ is large for the family $\cR$, the collection $\cP$ is nonempty.
By the Kuratowski--Zorn Lemma, there is a maximal element $\cR'$ in $\cP$.
We have $\Un\cR'\sub D$.
Assume that the set $S\sm \Un\cR'$ is in $p$.
Since $p$ is large for $\cR$, there is a set $R\in\cR$ such that $R\sub (S\sm\Un\cR')\cap A$.
Then $\cR'\subsetneq \cR'\cup\{R\}$, a contradiction.
Thus, $\Un\cR'\in p$.
\epf

For a set $A$ and a natural number $m$, let
\begin{align*}
 [A]^{\leq m}&:=\Un_{i\leq m}[A]^i,&
[A]^{< m}&:=\Un_{i<m}[A]^i.
\end{align*}

Let $S$ be an infinite semigroup and $n,m$ be natural numbers.
Let $\eseqint{a}{1}{n}\in S$ be a \emph{proper} sequence, i.e., $a_{H_1}\neq a_{H_2}$ for all sets $H_1,H_2\in\fin(\{1,\dotsc,n\})$ with $H_1<H_2$.
Define
\[
\SG^m(\eseqint{a}{1}{n}):=
\smallmedset{\{a_{H_1}, \dotsc, a_{H_m}\}}{\eseqint{H}{1}{m}\in\fin(\{1,\dotsc,n\})\text{ and }H_1<\dotsb <H_m}.
\]
If $m>n$, then the set $\SG^m(\eseqint{a}{1}{n})$ is empty.
For sets $\eseqint{F}{1}{n}\in \fin(S)$ such that all sequences in $\eprodint{F}{1}{n}$ are proper, define
\begin{gather*}
\SG^m[\eseqint{F}{1}{n}]:=\Un\smallmedset{\SG^m(\eseqint{a}{1}{n})}{(\eseqint{a}{1}{n})\in\eprodint{F}{1}{n}},\\
\SG^{\leq m}[\eseqint{F}{1}{n}]:=\Un_{i\leq m}\SG^i[\eseqint{F}{1}{n}],\quad \SG^{<m}[\eseqint{F}{1}{n}]:=\Un_{i< m}\SG^i[\eseqint{F}{1}{n}].
\end{gather*}

For nonempty sets $A,B$ in a semigroup $S$, let 
\[
A+B:=\sset{a+b}{a\in A, b\in B},
\]
and for an element $b\in S$, by $b+B$ we mean $\{b\}+B$.

\begin{proof}[Proof of Theorem~\ref{thm:mainfin}]
Fix natural numbers $m,k$ and let $\chi_m\colon [S]^m\to \{1,\dotsc,k\}$ be a coloring of  $[S]^m$.
Let $e\sub\cA$ be an idempotent in $\beta S$, large for the sequence $\eseq{\cR}$.
Since Theorem~\ref{thm:mainfin} for $m=2$, implies Theorem~\ref{thm:mainfin} for $m=1$, assume that $m\geq 2$. 
Define colorings $\chi_{m-1},\dotsc,\chi_1$ of the sets $[S]^{m-1},\dotsc,[S]^1$, respectively, as follows.
For each set $A\in[S]^{m-1}$, there is a unique number $\chi_{m-1}(A)\in\{1,\dotsc,k\}$ such that the set 
\[
\set{s\in S\sm A}{\chi_m(A\cup\{s\})=\chi_{m-1}(A)}
\]
is in $e$.
Analogously, define colorings $\chi_{m-2},\dotsc,\chi_1$ of the sets $[S]^{m-2},\dotsc,[S]^1$, respectively.

The function $\chi:=\Un_{i\leq m}\chi_m$ is a coloring of  $[S]^{\leq m}$. 
There is a $\chi$-monochromatic set $G$ in $e$ (strictly speaking, the set $G$ is $\chi_1$-monochromatic).
Assume that the color is green.
For each set $A\in [S]^{<m}$ which is green, the set $\sset{s\in S\sm A}{A\cup \{s\}\text{ is green}}$ is in $e$.
Thus, for each nonempty finite set $F\sub [S]^{<m}$ which is green, the set 
\[
G(F):=\bigcap_{A\in F}\set{s\in S\sm A}{A\cup \{s\}\text{ is green}}
\]
is in $e$, too.

Let $D\in e$.
Since $e$ is idempotent, the set 
\[
D^\star=
\sset{b\in D}{\Exists{B\in e}\End{B\sub D\wedge b+B\sub D}}
\]
is in $e$.

Define a strategy for \Alice{} as follows.

\textbf{1st round:} Let 
\[
D_1:=G.
\]
The set $D_1$ is in $e$.
By Lemma~\ref{lem:ult}, there is a family $\cR_1'\sub\cR_1$ of pairwise disjoint subsets of $D_1^\star$, whose union $A_1:=\Un\cR'_1$ is in $e$.  
\Alice{} plays the set $A_1$.
\Bob{} replies with a finite set $F_1\sub D_1^\star$.
Then there is a finite nonempty family $\cF_1\sub\cR_1'$ such that $F_1\sub \Un\cF_1$.
Let $V_1:=\Un\cF_1$.
We have $\SG^{\leq m}[V_1]=V_1\sub D_1$, and thus the set $\SG^{\leq m}[V_1]$ is green.

\textbf{2nd round:}
Since the set $V_1\sub D_1^\star$ is finite, there is a set $B_1\sub D_1$ in  $e$ such that $V_1+B_1\sub D_1$.
Let
\[
D_2:=B_1\cap G(\SG^{\leq m}[V_1]).
\]
Since  the set $\SG^{\leq m}[V_1]$ is green and finite, the set  $G(\SG^{\leq m}[V_1])$ is in $e$, and thus the set $D_2$ is in $e$, as well.
By Lemma~\ref{lem:ult}, there is a family $\cR_2'\sub\cR_2$ of pairwise disjoint subsets of $D_2^\star$, whose union $A_2:=\Un\cR'_2$ is in $e$.  
\Alice{} plays the set $A_2$ and 
\Bob{} replies with a finite set $F_2\sub A_2$.
Then there is a finite nonempty family $\cF_2\sub\cR_n'$ such that $F_2\sub\Un\cF_2$.
Let $V_2:=\Un\cF_2$.

The set $\SG^{\leq m}[V_1,V_2]$ is green:
By the previous step, the set $V_1$ is green, and the set $V_2$ is green as a subset of $D_2$.
Since $V_1+B_1\sub D_1$ and $V_2\sub B_1$, the set $V_1+V_2$ is green.
Thus, the set $\SG^1[V_1, V_2]=V_1\cup V_2\cup (V_1+V_2)$ is green.
Since $V_2\sub G(\SG^1[V_1])$, the set $\SG^2[V_1, V_2]$ is green.
Thus, the set $\SG^{\leq  m}[V_1, V_2]=\SG^1[V_1, V_2]\cup \SG^2[V_1, V_2]$ is green.

Fix a natural number $n\geq 2$.
Assume that the sets $\eseqint{D}{1}{n}$ in $e$ with $D_1\supseteq \dotsb\supseteq D_n$, finite nonempty families $\cF_1\sub\cR_1, \dotsc, \cF_n\sub\cR_n$ of pairwise disjoint sets, sets $A_1\sub D_1^\star,\dotsc, A_n\sub D_n^\star$ in $e$ and finite sets $F_1\sub A_1,\dotsc, F_n\sub A_n$ have already been defined such that for $V_1:=\Un\cF_1,\dotsc,V_n:=\Un\cF_n$, the set $\SG^{\leq m}[\eseqint{V}{1}{n}]$ is green.

\textbf{$\mathbf{(n+1)}$ round:}
Since the set $V_n\sub D_n^\star$ is finite, there is a set $B_n\sub D_n$ in $e$ such that $V_n+B_n\sub D_n$.
Let
\[
D_{n+1}:=B_n\cap G(\SG^{\leq m}[\eseqint{V}{1}{n}]).
\]
Since  the set $\SG^{\leq m}[\eseqint{V}{1}{n}]$ is green and finite, the set  $G(\SG^{\leq m}[\eseqint{V}{1}{n}])$ is in $e$.
Thus, the set $D_{n+1}$ is in $e$, too.
By Lemma~\ref{lem:ult}, there is a family $\cR_{n+1}'\sub\cR_{n+1}$ of pairwise disjoint subsets of $D_{n+1}^\star$, whose union $A_{n+1}:=\Un\cR'_{n+1}$ is in $e$.  
\Alice{} plays the set $A_{n+1}$ and 
\Bob{} replies with a finite set $F_{n+1}\sub A_{n+1}$.
Then there is a finite nonempty family $\cF_{n+1}\sub\cR'_{n+1}$ such that $F_{n+1}\sub\Un\cF_{n+1}$.
Let $V_{n+1}:=\Un\cF_{n+1}$.

\bclm\label{clm:main}
The set $\SG^{\leq m}[\eseqint{V}{1}{n+1}]$ is green.
\eclm

\bpf

Firstly, prove that for a natural number $j$, natural numbers $i_1<\dotsb<i_j\leq n+1$ and elements $a_{i_1}\in V_{i_1},\dotsc, a_{i_j}\in V_{i_j}$,
the element $b=a_{i_1}+\dotsb + a_{i_j}$ is in $D_{i_1}$:
If $j=1$, then $b= a_{i_1}\in V_{i_1}\sub D_{i_1}$.
Assume that the claim is true for a natural number $j$ and
$b=a_{i_1}+\dotsb + a_{i_{j+1}}$.
By the assumption, we have $b':=a_{i_2}+\dotsb + a_{i_{j+1}}\in D_{i_2}\sub D_{i_1+1}\sub B_{i_1}$.
Thus, $b=a_{i_1}+b' \in V_{i_1}+B_{i_1}\sub D_{i_1}$.

By the above, the set $\SG^1[\eseqint{V}{1}{n+1}]$ is green.
Let $A\in \SG^{\leq m}[\eseqint{V}{1}{n+1}]$ be a set such that $\card{A}>1$.
Then there are a natural number $l$ with $1<l\leq n+1$ and an element $b\in A$ such that $A\sm\{b\}\in\SG^{<m}[\eseqint{V}{1}{l}]$ and $b\in \SG^1[\eseqint{V}{l+1}{n+1}]$.
By the above, $b\in D_{l+1}\sub G(\SG^{<m}[\eseqint{V}{1}{l}])$, and thus the set $A$ is green.
\epf

(1) Since \Alice{} has no winning strategy in the game $\gfin(\cA,\cB)$, there is a play
\[
(A_1, F_1, A_2, F_2,\dotsc),
\]
where \Alice{} applies the strategy defined above, and the game is won by \Bob{}.
It provides that the set $\Un_{n\in\bbN}F_n$ is in $\cB$.

(2)
Let $\eseq{a}$ be a sequence in $\eprod{V}$.
Let $H_1,H_2\in\Fin$ be sets with $H_1< H_2$ and $n=\max H_1$.
Since $a_{H_2}\in \SG^1[V_{n+1},V_{n+2},\dotsc]$, we have $a_{H_2}\in D_{n+1}$.
Since $a_{H_1}\in\SG^1[\eseqint{F}{1}{n}]$ and $D_{n+1}\sub G(\SG^{<m}[\eseqint{H}{1}{n}])$, we have $a_{H_1}\neq a_{H_2}$. 

(3)
It follows from Claim~\ref{clm:main}.
\epf

\brem\label{rem:mainfin}
If $e$ is an idempotent from Theorem~\ref{thm:mainfin}
and $f\colon\fin(S)\to e$ is a function, then, in Theorem~\ref{thm:mainfin}, we may require that $\Un\cF_{n+1} \sub f(\Un_{i\leq n}\cF_i)$ for all $n$ or $R_{n+1}\sub f(\Un_{i\leq n}R_i)$ for all $n$, respectively.
\erem

\section{Applications}

\subsection{Natural numbers}
\label{ssec:nat}

Before we discuss consequences of Theorem~\ref{thm:mainfin}, for the semigroup $\bbN$ with the usual addition $+$, let us take a closer look to examples of sequences providing an idempotent in $\roth$.
Let $m,p,c\in\bbN$.
For each element $x\in\bbN^m$, let $S(m,p,c,x)$
be the set of all sums
\begin{align*}
cx_1 + \lambda_{2}x_{2} + \lambda_3x_3+\dotsb + \lambda_{m}x_{m},&\\
cx_2 + \lambda_{3}x_{3} + \dotsb + \lambda_{m}x_{m},&\\
\vdots&\\
cx_{m-1}+\lambda_mx_m,\\
cx_m,&
\end{align*}
where $|\lambda_2|,\dotsc, |\lambda_m|<p$.
An \emph{$(m,p,c)$-set} is the set $S(m,p,x)$ for some element $x\in\bbN^m$.
The importance of $(m,p,c)$-sets lies in the fact that they characterize partition regular systems of homogeneous linear equations~\cite[Section~3.3]{Rbook}. 

\bexm
\label{ex:mpc}
Let $\eseq{\cR}\sub\Fin$ be an enumeration of all families of $(m,p,c)$-sets.
By the result of Deuber and Hindman~\cite[Lemma~2]{DH} and the Ellis--Numakura Lemma, the family $\roth$ contains an idempotent, large for the sequence $\eseq{\cR}$.
\eexm

For other theorems in a similar spirit, we refer to the work of Bergelson and Hindman~\cite{BHcell}.

The following result is a consequence of Theorem~\ref{thm:mainfin}.
\bprp\label{prp:nat}
Assume that $m$ is a natural number and $\roth$ contains an idempotent, large for a sequence $\eseq{\cR}\sub\Fin$.
Then for each coloring of  $[\bbN]^m$, there are sets $R_1\in\cR_1, R_2\in\cR_2,\dotsc$ such that
$R_1<R_2<\dotsb$ and the partite $m$-sumgraph of the sequence $\eseq{R}$ is monochromatic.
\eprp

\bpf
Bob has a winning strategy in the game $\gone(\roth,\roth)$.
Apply Theorem~\ref{thm:mainfin} to the families $\cA=\cB=\roth$ and to the sequence $\eseq{\cR}$.
Then there are sets $R_1\in\cR_1, R_2\in\cR_2,\dotsc$ such that all sequences in the product $\eprod{R}$ are proper and the partite $m$-sumgraph of the sequence $\eseq{R}$ is monochromatic.
For each finite subset of $\bbN$, there is a natural number $n$ such that the set is disjoint with $\{n,n+1,\dotsc\}$, an element of the idempotent.
By Remark~\ref{rem:mainfin}, we may assume that $R_1<R_2<\dotsb$.
\epf

By Example~\ref{ex:mpc}, the Deuber--Hindman Theorem~\ref{thm:DH} follows from Proposition~\ref{prp:nat}.
In our terminology, the claim of Theorem~\ref{thm:DH} can be replaced by: the $1$-partite sumgraph of the sequence $\eseq{R}$ is monochromatic.
In the original formulation of Theorem~\ref{thm:DH}, the sequence $\eseq{\cR}$ is an enumeration of all families of $(m,p,c)$-sets.

By Example~\ref{ex:mpc}, the family $\roth$ contains an idempotent for the sequence $\eseq{\cR}$ of families of all arithmetic progressions of length $1,2,\dotsc$, respectively.
Thus, we have the following corollary from Proposition~\ref{prp:nat}.

\bcor
Assume that $m$ is a natural number.
Then for each coloring of  $[\bbN]^m$, there are arithmetic progressions $\eseq{R}$ of length $1,2\dotsc$, respectively, such that all sequences in $\eprod{R}$ are proper 
and the partite $m$-sumgraph of the sequence $\eseq{R}$ is monochromatic.
\ecor

Also, mentioned in the introduction Theorem~\ref{thm:BHweak} and its higher-dimensional version, below Theorem~\ref{thm:BHweakHD}, follows from Proposition~\ref{prp:nat}.

\bthm[{Bergelson, Hindman~\cite[Theorem~2.5.]{BH}}]
\label{thm:BHweakHD}
Assume that $m$ is a natural number and $\roth$ contains an idempotent, large for a sequence $\eseq{\cR}\sub\Fin$.
Then for each coloring of  $[\bbN]^m$, there are increasing sequences $x^{(n)}_1,x^{(n)}_2,\dotsc$ of natural numbers and pairwise disjoint sets $R_1\in\cR_1, R_2\in\cR_2,\dotsc$ with the following property.
For every natural numbers $n_1<\dotsb<n_m$ and elements
\begin{align*}
a_1&\in R_{n_1}\cup\Un\sset{\FS(x^{(l)}_{1},\dotsc,x^{(l)}_{n_1})}{l\leq n_1},\\
&\vdots\\
a_m&\in R_{n_m}\cup\Un\sset{\FS(x^{(l)}_{n_{m-1}+1},\dotsc,x^{(l)}_{n_m})}{n_{m-1}<l\leq n_m},
\end{align*}
with $a_1<\dotsb<a_m$, the sets $\{a_1,\dotsc,a_m\}$ have the same color.
\ethm

\bpf
Fix a coloring of  $[\bbN]^m$.
Let $\eseq{R}$ be sets from Proposition~\ref{prp:nat}
and $\eseq{x}$ be any sequence in $\eprod{R}$.
For every natural numbers $n_1<\dotsb<n_m$ and elements
\begin{align*}
a_1&\in R_{n_1}\cup\Un\sset{\FS(x_{1},\dotsc,x_{n_1})}{l\leq n_1},\\
&\vdots\\
a_m&\in R_{n_m}\cup\Un\sset{\FS(x_{n_{m-1}+1},\dotsc,x_{n_m})}{n_{m-1}<l\leq n_m},
\end{align*}
the set $\{a_1,\dotsc,a_m\}$ belongs to the partite $m$-sumgraph of the sequence $\eseq{R}$, a monochromatic set.
\epf

Theorem~\ref{thm:BHweakHD}, in the original form,  is discussed in Subsection~\ref{ssec:BH}.

\subsection{Covers}
\label{ssec:covers}

Let $\cA,\cB$ be families of sets.
Define
\smallskip

\begin{labeling}{$\sone(\cA,\cB)$:}
\item [$\sone(\cA,\cB)$:]  for each sequence $\eseq{A}\in\cA$, there are elements $a_1\in A_1, a_2\in A_2,\ldots$ such that $\sset{a_n}{n\in\bbN}\in \cB$.
\end{labeling}
\smallskip

Some topological covering properties can be expressed, using the above pattern.
A space is \emph{Rothberger} if it satisfies $\soo{X}$~\cite{rothb}.
For sets of reals, the Rothberger property implies \emph{strong measure zero} and it can be viewed as a topological version of strong measure zero.
A counterpart of Theorem~\ref{thm:hur}, for the Rothberger property, is the famous result of Pawlikowski.

\bthm[{Pawlikowski~\cite{Paw94}}]\label{thm:paw}
A space $X$ satisfies the Rothberger property $\soo{X}$ if and only if \Alice{} has no winning strategy in the game $\goo{X}$. 
\ethm

For a separable metrizable space, \Bob{} has a winning strategy in the game $\goo{X}$ if and only if the space is countable~\cite{Paw94}.
By the results of Sierpi\'{n}ski~\cite{Hure27} and Laver~\cite{Laver}, existence of an uncountable Rothberger set of reals is independent from ZFC.
Thus, in the realm of sets of reals, indeterminacy of the game $\goo{X}$ is independent from ZFC.
As the Menger property, also the Rothberger property has a connection with local properties in function spaces.
A space $Y$ has \emph{countable strong fan tightness} if for every point $y\in Y$ and a sequence $\eseq{A}\sub Y$ such that $y\in\bigcap_{n\in\bbN}\overline{A_n}$, there are elements $a_1\in A_1, a_2\in A_2,\dotsc$ such that $y\in\overline{\sset{a_n}{n\in\bbN}}$.
Let $X$ be a set of reals.
By the result of Sakai~\cite{sakai}, the space $\Cp(X)$ has countable strong fan tightness if and only if the space $X$ is Rothberger in all finite powers, equivalently, it satisfies $\sww{X}$.
Another local properties in function spaces were considered by Gerlits and Nagy~\cite{gn}.
An infinite cover of a space is a \emph{$\gamma$-cover}, if each element of the space belongs to all but finitely many sets in the cover.
For a space, let $\Ga$ be the family of all open $\gamma$-covers of the space.
Gerlits and Nagy proved that a space $X$ satisfies $\swg{X}$ (also known as the property $\gamma$) if and only if the space $\Cp(X)$ has the Fr\'{e}chet--Urysohn property, i.e., each point in the closure of a set in $\Cp(X)$ is a limit of a sequence of points from that set.

Diagram~\ref{diag} (a part of so--called the Scheepers diagram~\cite{Wiki}) shows the relations between considered properties.
In order to study the above properties in a coherent manner, we replace the Rothberger property $\soo{X}$ by its equivalent property $\swl{X}$~\cite[Theorem~17]{coc1}, and the Menger property $\sfinoo{X}$ by $\sfinwlx{X}$~\cite[Figure~2]{coc2}.

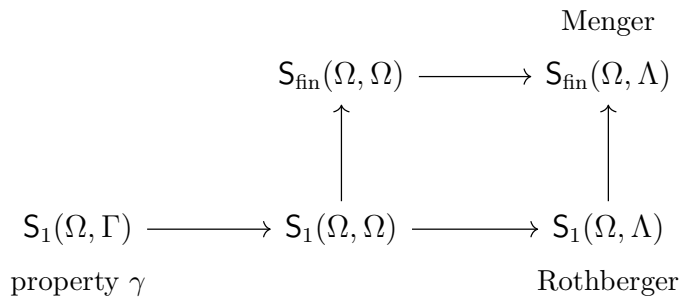
\begin{figure}[h]
		\begin{tikzcd}[ampersand replacement=\&,column sep=.7cm, row sep=1.3cm]
			{}\&\&\&\&\text{\small{Menger}}\\[-3.1em]
			{}\&\&\sfinww{X}\arrow{rr}\&\&\sfinwlx{X}\\
			{\swg{X}}\arrow{rr}\&\&\sww{X}\arrow{rr}\arrow{u}\&\&{\swl{X}}\arrow{u}\\ [-3.1em]
			{\text{\small{property $\gamma$}}}\&\&\&\&{\text{\small{Rothberger}}}
		\end{tikzcd}

\caption{Relations between considered combinatorial covering properties}
\label{diag}
\end{figure}

Using Theorems~\ref{thm:hur} and~\ref{thm:paw}, Scheepers showed that considered covering properties can be characterized in the language of games and winning strategies (\cite[Theorems~3 and~5]{schpart}, \cite[Theorems~13 and~35]{coc6}, \cite[Theorem~30]{coc1}).

\bthm[{Scheepers~\cite{schpart, coc6, coc1}}]\label{thm:gameequiv}
Let $X$ be a separable metrizable space and $\cB\in\{\Op, \Lambda, \Om,\Ga\}$.
\be 
\item 
The space $X$ satisfies $\sfin(\om{X},\cB)$ if and only if \Alice{} has no winning strategy in the game $\gfin(\om{X},\cB)$. 
\item The space $X$ satisfies $\sone(\om{X},\cB)$ if and only if \Alice{} has no winning strategy in the game $\gone(\om{X},\cB)$. 
\ee

\ethm

Theorem~\ref{thm:gameequiv} is a bridge between properties and games.
Since we want to apply Theorem~\ref{thm:mainfin}, in the context of open covers of spaces, we need some additional preparation.
For a space $X$, we introduce a semigroup operation on $\tau$ and give examples of families $\cA\sub\tau$ containing an idempotent, large for a sequence $\eseq{\cR}\sub\fin(\tau)$.
To this end, we use observations of Tsaban~\cite{tsramsey}.

\bdfn
Let $S$ be an infinite semigroup.
The family $\cA\sub\rothx{S}$ is a \emph{superfilter} on $S$ if
\be
\item
for each set $A \in\cA$, all subsets of $S$ that contain $A$ are in $\cA$,
\item whenever $A_1 \cup A_2\in \cA$, then $A_1$ or $A_2$ is in A; equivalently, for each
set $A \in\cA$ and each coloring of $A$, there is in $\cA$ a monochromatic subset of $A$.
\ee
\edfn

\bexm
For a space $X$, the family $\om{X}$ is a superfilter on $\tau\sm\{X\}$.
\eexm

\bdfn[{\cite[Definition2.4.]{tsramsey}}]
A superfilter $\cA$ on a semigroup $S$ is an \emph{idempotent superfilter} if for each set $A\sub S$, if the set 
\[
\smallmedset{b\in S}{\Exists{C\in\cA}\End{b+C\sub A}}
\]
is in $\cA$, then $A\in\cA$. 
\edfn

\bexm
\label{exm:idem}
Let $S$ be an infinite semigroup.
\be

\item 
By the Ellis--Numakura Lemma there is an idempotent ultrafilter in $\beta S$, an idempotent superfilter.

\item
\label{itm:trinv}
Each \emph{translation-invariant} superfilter $\cA$ on $S$, that is a superfilter such that for all $s\in S$ and $A\in\cA$, we have $s+A\in\cA$, is idempotent. 
\item \label{itm:max}If $S=\sset{a_n}{n\in\bbN}$, where elements $\eseq{a}$ are pairwise different and it has the semigroup operation $\vee$ defined by $a_n\vee a_m: = a_{\max(n,m)}$, then every superfilter on $S$ is translation-invariant, and thus idempotent.
\item If $S=\sset{a_n}{n\in\bbN}$, where elements $\eseq{a}$ are pairwise different and it has the semigroup operation $\wedge$ defined by $a_n\wedge a_m: = a_{\min(n,m)}$, then  $\rothx{S}$ is an idempotent superfilter on $S$ but it is not translation-invariant.



\item
Let $X$ be a space, $\cU\sub\tau$ and $\FS(\cU)$ be the family of all finite unions of elements from $\cU$.
Assume that $\cU$ is an open cover of $X$ with no finite subcover and $S:=\FS(\cU)$ is a subsemigroup of $(\tau,\cup)$.
Then the family 
\[
\cA:=\sset{A\sub S}{A\in \om{X}}
\]
is an idempotent superfilter on $S$.
\ee
\eexm

Let $\eseq{a}$ be a sequence in a semigroup $S$. For a natural number $n$, let
\[
\FS(a_n,a_{n+1},\dotsc):=\sset{a_H}{H\in\fin(\{n,n+1,\dotsc\})}.
\]
Identify a set $A$ with $[A]^1$.

\bexm
\label{exm:tsidem}
Let $\eseq{a}$ be a proper sequence in a semigroup $S$ and $\cA$ be an idempotent superfilter on $S$ such that $\FS(\eseqstart{a}{n})\in\cA$ for all $n$.
Define $\cR_n:=\FS(a_n,a_{n+1},\dotsc)$ for all $n$.
By the results of Tsaban~\cite[Lemma~3.4., Theorem~2.7.]{tsramsey}, the superfilter $\cA$ contains an idempotent, large for the sequence $\eseq{\cR}$.
\eexm

We derive the following results from Theorem~\ref{thm:mainfin}.

\bprp\label{prp:corfin}
Let $S$ be a countable set and $m$ be a natural number.
Assume that $\cA$ is a superfilter on $S$ and $\cB$ is a family of subsets of $S$, closed under taking supersets, such that \Alice{} has no winning strategy in the game $\gfin(\cA,\cB)$.
Then for each coloring of $[S]^m$, there are pairwise disjoint sets $\eseq{F}\in\fin(S)$ with the following properties.
\be
\item The set $\Un_{n\in\bbN}F_n$ is in $\cB$.
\item The partite $m$-graph of the sequence $\eseq{F}$ is monochromatic.
\ee
\eprp

\bpf
Let $S=\sset{s_n}{n\in\bbN}$, where elements $\eseq{s}$ are pairwise different.
Consider $S$ as a semigroup with the semigroup operation $\vee$ defined as in Example~\ref{exm:idem}(\ref{itm:max}).
By Example~\ref{exm:idem}(\ref{itm:max}), the family $\cA$ is an idempotent superfilter on $S$.
Let $\cR_n:=\{a_n,a_{n+1},\dotsc\}$ for all $n$.
By the definition of $\vee$, we have $\cR_n=\FS(a_n, a_{n+1},\dotsc)\in\cA$ for all $n$.
By Example~\ref{exm:tsidem}, the superfilter $\cA$ contains an idempotent, large for the sequence $\eseq{\cR}$.
Fix a coloring of $[S]^m$.
Let $\eseq{F}\in\fin(S)$ be sets, obtained from Theorem~\ref{thm:mainfin}.
Since all sequences in $\eprod{F}$ are proper, the sets $\eseq{F}$ are pairwise disjoint.
The set $\Un_{n\in\bbN}F_n$ is in $\cB$ and the partite $m$-graph of the sequence $\eseq{F}$ is equal to the $m$-partite sumgraph of that sequence, a monochromatic set.
\epf

\bprp\label{prp:corone}
Let $S$ be a countable set and $m$ be a natural number.
Assume that $\cA$ is a superfilter on $S$ and $\cB$ is a family of subsets of $S$ such that  \Alice{} has no winning strategy in the game $\gone(\cA,\cB)$.
Then for each coloring of $[S]^m$, there are pairwise different elements $\eseq{a}\in S$ with the following properties.
\be
\item The set $\sset{a_n}{n\in\bbN}$ is in $\cB$.
\item The $m$-graph $[\sset{a_n}{n\in\bbN}]^m$ of the sequence $\eseq{a}$ is monochromatic.
\ee
\eprp

\bpf
The proof is similar to the proof of Proposition~\ref{prp:corfin}.
\epf

\bdfn
A family $\cB$ of open covers of a space is  \emph{semiregular} if for every family $\cU\in \cB$ and a finite-to-one function $f\colon\cU\to\tau$
with $U \sub f(U)$ for all $U\in\cU$, the image of $f$ is in $\cB$.
\edfn

\bexm
For a space, the families $\opn{X}, \Lambda, \om{X}$ are semiregular.
\eexm

\blem\label{lem:color}
Let $X$ be a space and $\cB$ be a semiregular family of covers of $X$.
Assume that for every family $\cU\in\om{X}$ and coloring of  $[\tau]^2$, there are pairwise disjoint sets $\eseq{\cF}\in\fin(\cU)$ such that the family $\Un_{n\in\bbN}\cF_n$ is in $\cB$ and the partite graph of the sequence $\eseq{\cF}$ is monochromatic.
Then the space  satisfies $\sfin(\om{X},\cB)$.
Moreover, if the sets $\eseq{\cF}$ are singletons, then the space satisfies $\sone(\om{X},\cB)$.
\elem

\bpf
Let $\eseq{\cU}\in\om{X}$.
Assume that the family  $\cU':=\bigcap_{n\in\bbN}\cU_n$ is in $\om{X}$.
By the assumption, there are pairwise disjoint sets $\eseq{\cF}\in\fin(\cU')$ such that the family $\Un_{n\in\bbN}\cF_n$ is in $\cB$.
We have also $\cF_1\sub\cU_1, \cF_2\sub\cU_2,\dotsc$.
Thus, the space $X$ satisfies $\sfin(\om{X},\cB)$.

Now assume that the above family $\cU'$ is not in $\om{X}$. 
Since $\om{X}$ is a superfilter, the families $\cU_n\sm \cU'$ are in $\om{X}$ for all $n$.
Thus, we may assume that the family $\cU'$ is empty.
Let $\sset{W_n}{n\in\bbN}$ be an increasing open cover of $X$.
Then the family
\[
\cU:=\Un_{n\in\bbN} \sset{U\cap W_n}{U\in \cU_n}
\]
is in $\om{X}$.
For each $\{U,V\}\in [\tau]^2$, let 
\[
\chi(\{U,V\}) :=
\begin{cases}
1,&\text{ if }U,V\in \sset{U\cap W_n}{U\in \cU_n}\text{ for some }n,\\
2, &\text{ otherwise.}
\end{cases}
\]
By the assumption, there are pairwise disjoint sets $\eseq{\cF}\in\fin(\cU)$ such that the family $\Un_{k\in\bbN}\cF_k$ is in $\cB$ 
and the partite graph of the sequence $\eseq{\cF}$ is monochromatic.

Assume that the color is $1$.
Then there is a natural number $n$ such that $\Un_{k\in\bbN}\cF_k\sub \sset{U\cap W_n}{U\in \cU_n}$, and thus the family $\Un_{k\in\bbN}\cF_k$ is not a cover of $X$, a contradiction.
We conclude that the color is $2$.
Fix natural numbers $i,j$ with $i\neq j$ and sets $U\in\cF_i$, $V\in\cF_j$.
Since $\chi(\{U,V\})=2$, the sets $U$ and $V$ comes from different families $\sset{U\cap W_n}{U\in \cU_n}$.
It implies that the families 
\[
\Un_{k\in\bbN}\cF_k\cap \sset{U\cap W_n}{U\in \cU_n}
\]
are finite for all $n$.
Then there are finite sets $\tilde{\cF}_n\sub\cU_n$ for all natural numbers $n$ and a finite-to-one surjection
\[
f\colon \Un_{k\in\bbN}\tilde{\cF}_k \to \Un_{n\in\bbN}\cF_n.
\]
Since the family $\cB$ is semiregular, the family $\Un_{n\in\bbN}\cF_n$ is in $\cB$.
Thus, the space satisfies $\sfin(\om{X},\cB)$.
\epf

The following result gives characterizations of covering properties using colorings.

\bthm\label{thm:corfin}
Let $X$ be a separable metrizable space and $\cB\in \{\Op,\Lambda, \Om\}$.
The following assertions are equivalent.
 
\be
\item The space $X$ satisfies $\sfin(\om{X},\cB)$.
\item \Alice{} has no winning strategy in the game $\gfin(\om{X}, \cB)$.
\item For every natural number $m$, coloring of  $[\tau]^m$ and family $\cU\in\om{X}$, there are pairwise disjoint sets $\eseq{\cF}\in\fin(\cU)$ such that the family $\Un_{n\in\bbN}\cF_n$ is in  $\cB$ and the partite $m$-graph of the sequence $\eseq{\cF}$ is monochromatic.
\item For every coloring of  $[\tau]^2$ and family $\cU\in\om{X}$, there are pairwise disjoint sets $\eseq{\cF}\in\fin(\cU)$ such that the family $\Un_{n\in\bbN}\cF_n$ is in  $\cB$ and the partite graph of the sequence $\eseq{\cF}$ is monochromatic.
\ee
\ethm

\bpf
(1)$\Impl$(2) Apply Theorem~\ref{thm:gameequiv}(1).

(2)$\Impl$(3) Apply Proposition~\ref{prp:corfin}.

(3)$\Impl$(4) It is straightforward.

(4)$\Impl$(1) Apply Lemma~\ref{lem:color}.
\epf

For $\cB=\Om$, the implication (1)$\Impl$(4), in Theorem~\ref{thm:corfin}, is the result of Scheepers~\cite[Theorem~10]{coc1} and the implication (4)$\Impl$(1) is the result of Just, Miller, Scheepers and Szeptycki \cite[Theorem~6.2]{coc2}.
For $\cB=\Lambda$, the equivalence (1)$\Leftrightarrow$(4) is the result of Scheepers~\cite[Theorem~6]{schpart}.

\bthm
\label{thm:corone}
Let $X$ be a separable metrizable space and $\cB\in \{\Lambda, \Om, \Ga\}$.
The following assertions are equivalent.
\be
\item The space $X$ satisfies $\sone(\om{X},\cB)$.
\item \Alice{} has no winning strategy in the game $\gone(\om{X}, \cB)$.
\item For every natural number $m$, coloring of  $[\tau]^m$ and family $\cU\in\om{X}$, there are pairwise different sets $\eseq{U}\in\cU$ such that the family $\sset{U_n}{n\in\bbN}$ is in $\cB$ and the $m$-graph $[\sset{U_n}{n\in\bbN}]^m$ of the sequence $\eseq{U}$ is monochromatic.
\item For every coloring of  $[\tau]^2$ and family $\cU\in\om{X}$, there are pairwise different sets $\eseq{U}\in\cU$ such that the family $\sset{U_n}{n\in\bbN}$ is in $\cB$ and the graph $[\sset{U_n}{n\in\bbN}]^2$ is monochromatic.
\ee
\ethm

\bpf
(1)$\Impl$(2) Apply Theorem~\ref{thm:gameequiv}(2).

(2)$\Impl$(3) Apply Proposition~\ref{prp:corfin}.

(3)$\Impl$(4) It is straightforward.

(4)$\Impl$(1) Apply Lemma~\ref{lem:color}.
\epf

For $\cB=\Lambda$, the equivalence (1)$\Leftrightarrow$(4) is the result of Scheepers~\cite[Theorem~4]{schpart}.
For $\cB\in\{\Om,\Ga\}$, Theorem~\ref{thm:corone} captures another results of Scheepers~\cite[Theorems~24, 25 and~30]{coc1}.

\subsection{Local properties}
\label{ssec:loc}

Let $Y$ be a space and $y\in Y$.
Define 
\[
\Om_y:=\sset{A\sub Y}{y\in\overline{A}}\text{ and }\Ga_y:=\sset{A\sub Y}{\Exists{\eseq{a}\in A}\End{\lim_{n\to\infty}a_n=y}}.
\]

The family $\Om_y$ is a superfilter on $\Pof(Y)$.

Let $X$ be a space and $\bfO\in\Cp(X)$ be the constant zero function.
The space $\Cp(X)$ is 

\bi
\item strongly Fr\'{e}chet--Urysohn if it satisfies $\sone(\Om_\bfO,\Ga_\bfO)$ (equivalently it is Fr\'{e}chet--Urysohn),
\item has countable strong fan tightness if it satisfies $\sone(\Om_\bfO,\Om_\bfO)$,
\item has countable strong tightness if it satisfies $\sone(\Om_\bfO,\Om_\bfO)$.
\ei

Using similar techniques as in Subsection~\ref{ssec:covers} we derive the following results.

\bthm\label{thm:corfinloc}
Let $X$ be a set of reals.
The following assertions are equivalent.
 
\be
\item The space $\Cp(X)$ satisfies $\sfin(\Om_\bfO, \Om_\bfO)$.
\item \Alice{} has no winning strategy in the game $\gfin(\Om_\bfO, \Om_\bfO)$.
\item For every natural number $m$, coloring of  $[\Pof(Y)]^m$ and set $A\in\Om_\bfO$, there are pairwise disjoint sets $\eseq{F}\in\fin(A)$ such that the set $\Un_{n\in\bbN}F_n$ is in $\Om_\bfO$ and the partite $m$-graph of the sequence $\eseq{F}$ is monochromatic.
\item For every coloring of  $[\Pof(Y)]^2$ and set $A\in\Om_\bfO$, there are pairwise disjoint sets $\eseq{F}\in\fin(A)$ such that the set $\Un_{n\in\bbN}F_n$ is in $\Om_\bfO$ and the partite graph of the sequence $\eseq{F}$ is monochromatic.
\ee
\ethm

\bthm\label{thm:coroneloc}
Let $X$ be a set of reals. 
For the space $\Cp(X)$, the following assertions are equivalent.
For the space $\Cp(X)$ and $\cB\in \{\Om_\bfO, \Ga_\bfO\}$ the following assertions are equivalent.
 
\be
\item The space $Y$ satisfies $\sone(\Om_\bfO, \cB)$
\item \Alice{} has no winning strategy in the game $\gone(\Om_\bfO, \cB)$.
\item For every natural number $m$, coloring of  $[\Pof(Y)]^m$ and set $A\in\cB$, there are pairwise different elements $\eseq{a}\in A$ such that the set $\sset{a_n}{n\in\bbN}$ is in $\cB$ and the $m$-graph $[\sset{a_n}{n\in\bbN}]^m$ of the sequence $\eseq{a}$ is monochromatic.
\item For every coloring of  $[\Pof(Y)]^2$ and set $A\in\Om_\bfO$, there are pairwise different elements $\eseq{a}\in A$ such that the set $\sset{a_n}{n\in\bbN}$ is in $\cB$ and the graph $[\sset{a_n}{n\in\bbN}]^2$ of the sequence $\eseq{a}$ is monochromatic.
\ee
\ethm

\section{Richer monochromatic structures}

\subsection{Around the Tsaban and Scheepers Theorems}
\label{ssec:ts}

Let $X$ be a separable Menger space and $\cup$ be the semigroup operation on $\tau$.
Fix a coloring of $[\tau]^2$ and a family $\cU\in\om{X}$.
Comparing Theorems~\ref{thm:ts} and~\ref{thm:corfin}, the question arises, whether there are pairwise disjoint sets $\eseq{\cF}\in\fin(\cU)$ such that the family $\Un_{n\in\bbN}\cF_n$ forms a cover of $X$ and both, the sumgraph of the sequence $\eseq{\Un\cF}$ and the partite graph of the sequence $\eseq{\cF}$ have the same color.
Another question is whether, in Theorem~\ref{thm:corfin}(4), we can request that the partite sumgraph of the sequence $\eseq{\cF}$ is monochromatic.
The following example shows that, even for countable discrete spaces, the answers are negative.

\bexm
Let $\Fin$ be a discrete space.
For each natural number $n$, let $U_n:=\sset{F\in\Fin}{n\notin F}$.
Fix a finite family $\cF\sub\Fin$. There is a natural number $n$ such that $n\notin\Un\cF$, and thus $\cF\sub U_n$.
Then the family $\cU:=\sset{U_n}{n\in\bbN}$ is in $\Om$.
Consider the subsemigroup $\FS(\cU)$ of the semigroup $(\tau,\cup)$.
Define a coloring $\chi\colon[\FS(\cU)]^2\to \{0,1\}$ by
\[
\chi(\{U,V\}):=
\begin{cases}
0,&\text{ if }U,V\in\cU,\\
1,&\text{ otherwise.}
\end{cases}
\]
For natural numbers $n,m$ with $n\neq m$, we have $U_n\cup U_m\notin \cU$.
Consequently, for each sequence $n_1<n_2<\dotsb$ of natural numbers,
the sequence $U_{n_1},U_{n_2},\dotsc$ is proper and the sum graph of that sequence is not monochromatic.
\eexm

In this subsection, we combine Theorem~\ref{thm:mainfin} with Theorems~\ref{thm:corfin},~\ref{thm:corone}.

A cover of a space is \emph{ascending} if it contains a cover  $\sset{U_n}{n\in\bbN}$ with $U_1\subsetneq U_2\subsetneq\dotsb$.
Let $\Asc$ be the family of all open ascending covers of a space.

\bthm\label{thm:sfinww}
Let $X$ be a separable metrizable space and $\cB\in\{\Lambda,\Om\}$.
Assume that $X$ satisfies $\sfin(\Om,\cB)$ and $\cup$ is the semigroup operation on $\tau$.
Then for every coloring of $[\tau]^2$ and family $\cU\in\Om$ with no finite subcover, there are pairwise disjoint sets $\eseq{\cF}\in\fin(\cU)$ with the following properties.
\be
\item The family $\Un_{n\in\bbN}\cF_n$ is in $\cB$.
\item The sequence $\eseq{\Un\cF}$ is proper.
\item 
The sumgraph of the sequence $\eseq{\Un\cF}$ is monochromatic.
\ee
\ethm

\begin{proof}[{Proof of Theorem~\ref{thm:sfinww}}]
Fix a coloring $\chi$ of $[\tau]^2$.
By Theorem~\ref{thm:corfin}(4), there are pairwise disjoint families $\eseq{\cH}\in\fin(\cU)$ such that the family $\Un_{n\in\bbN}\cH_n$  is in $\cB$ and the partite graph of the sequence $\eseq{\cH}$ is $\chi$-monochromatic.

Let $S:=\fin(\Un_{n\in\bbN}\cH_n)$ be a subsemigroup of $\tau$ and $\cA$ be the family of all sets $A\sub S$ containing a sequence of sets $\cF_1\subsetneq \cF_2\subsetneq\dotsb$ such that the family $\Un_{n\in\bbN}\cF_n$ is in $\cB$.
The family $\cA$ is an idempotent superfilter on $S$:
Let $A\cup B\in \cA$.
Then $A\cup B$ contains a sequence of sets $\cF_1\subsetneq \cF_2\subsetneq\dotsb$ such that the family $\Un_{n\in\bbN}\cF_n$ is in $\cB$, and thus, $A$ or $B$, contains a subsequence of $\eseq{\cF}$.
Then $A\in\cA$ or $B\in\cA$.
The superfilter $\cA$ is translation-invariant.
By Example~\ref{exm:idem}(\ref{itm:trinv}), the superfilter $\cA$ is idempotent.

Since the family $\cU$ has no finite subcover, for each natural number $n$, there is an element $x_n\in X\sm\Un_{i\leq n}\Un\cH_i$.
For each $n$, let
\[
\cV_n:=\smallmedset{ \cF\in\fin(\Un_{i\geq n}\cH_i) }{ \eseqint{x}{1}{n-1}\in\Un\cF}.
\] 
We have $\bigcap_{n\in\bbN}\cV_n=\emptyset$.
By the result of Tsaban~\cite[Lemma~3.4., Theorem~2.7.]{tsramsey}, the superfilter $\cA$ contains an idempotent, large for the sequence $\eseq{\cV}$.
Let $\cB'$ be the family of all sets $A\sub S$ whose union is in $\cB$.
For each set $A\in\cA$, the family $\sset{\Un\cF}{\cF\in A}$ is in $\Asc$.
Since the space $X$ satisfies $\sfin(\Om,\cB)$, it also satisfies $\sone(\Asc, \cB)$, and thus, \Alice{} has no winning strategy in the game $\gone(\Asc,\cB)$~\cite[Proposition~4.5., Theorem~5.5.]{tsramsey}.
Consequently, \Alice{} has no winning strategy in the game $\gone(\cA, \cB')$.

Define a coloring $\kappa$ of $\fin(\cU)$ such that $\kappa(\cF):=\chi(\Un\cF)$ for all $\cF\in\fin(\cU)$.
By Theorem~\ref{thm:mainfin}, there are sets $\eseq{\cF}\in\fin(\cU)$ such that $\sset{\cF_n}{n\in\bbN}$ is in $\cB$, the sequence $\eseq{\cF}$ is proper in $S$ and the sumgraph of the sequence $\eseq{\cF}$ is $\kappa$-monochromatic.
For each finite subset of $S$, there is a natural number $n$ such that the set is disjoint with $\cV_n$, an element of the idempotent.
By Remark~\ref{rem:mainfin}, we may assume that there is a sequence $1=i_1<i_2<\dotsb$ of natural numbers such that
\[
\cF_n\in \cV_{i_n}\sm \cV_{i_{n+1}}
\] 
for all $n$.
Then the partite graph of the sequence $\eseq{\cF}$ is $\chi$-monochromatic.
By the definition of sets $\cV_n$ the sequence $\eseq{\Un\cF}$ is proper in $\tau$ and the sumgraph of the sequence $\eseq{\Un\cF}$ is $\chi$-monochromatic.
\epf

\subsection{Around the Bergelson--Hindman Theorem}
\label{ssec:BH}


Let $\cdot$ be the usual multiplication on $\bbN$.
Let $\eseq{y}$ be a sequence of natural numbers.
For a set $H=\{i_1,\dotsc,i_n\}\in\Fin$ with $i_1<\dotsb<i_n$, where $n$ is a natural number, define
\[
y^\odot_H:=y_{i_1} \dotsm  y_{i_n}.
\]

Bergelson and Hindman, proved the following beautiful theorem engaging both, addition and multiplication on $\bbN$.
The forthcoming Theorem~\ref{thm:BHvSz} is a generalization of this result.

\bthm[{Bergelson, Hindman~\cite[Theorem~2.5.]{BH}}]
\label{thm:BHfullHD}
Assume that $m$ is a natural number and $\roth$ contains an idempotent, large for a sequence $\eseq{\cR}\sub\Fin$.
Then for each coloring of  $[\bbN]^m$, there are increasing sequences $x^{(n)}_1,x^{(n)}_2,\dotsc$ for all $n$, and $\eseq{y}$ of natural numbers and pairwise disjoint sets $R_1\in\cR_1, R_2\in\cR_2,\dotsc$ with the following property.
For every natural numbers $n_1<\dotsb<n_m$ and elements
\begin{align*}
a_1&\in R_{n_1}\cup\Un\sset{\FS(x^{(l)}_{1},\dotsc,x^{(l)}_{n_1})}{l\leq n_1}\cup\sset{y^\odot_H}{H\in\fin(\{1,\dotsc, n_1\})},\\
&\vdots\\
a_m&\in R_{n_m}\cup\Un\sset{\FS(x^{(l)}_{n_{m-1}+1},\dotsc,x^{(l)}_{n_m})}{n_{m-1}<l\leq n_m}\cup\sset{y^\odot_H}{H\in\fin(\{n_{m-1},\dotsc, n_m\})},
\end{align*}
with $a_1<\dotsb<a_m$, the sets $\{a_1,\dotsc,a_m\}$ have the same color.
\ethm

An ultrafilter on $\bbN$ is \emph{multiplicative} if for each set $A$ in the ultrafilter, there is a natural number $y$ such that the set $A/y:=\sset{x\in\bbN}{x\cdot y\in A}$ is in the ultrafilter.
Let $\eseq{F}\in \Fin,$ $\eseq{y}$ be a sequence of natural numbers and $H\in\Fin$.
Define
\[
F(H):=
\begin{cases}
F_{i_1},& \text{if }n=1,\\
\sset{y^\odot_{H'}}{H'\in\fin(H)},& \text{otherwise}.
\end{cases}
\]
Fix natural numbers $n,m$.
Let $\eseqint{F}{1}{n}\in\Fin$ be a sequence such that all sequences in $F(H_1)\x\dotsb\x F(H_k)$, where $k$ is a natural number, $\eseqint{H}{1}{k}\in\Fin$ and $H_1<\dotsb<H_k$, are proper (with respect to $+$).
Let $\tSG^m[\eseqint{F}{1}{n}]$ be the union of all sumgraphs $\SG^m(F(H_1),\dotsc,F(H_k))$, where $k$ is a natural number,  $H_1,\dotsc,H_m\in\fin(\{1\dotsc,k\})$ and $H_1<\dotsb<H_k$.
Define
\begin{gather*}
\tSG^{\leq m}[\eseqint{F}{1}{n}]:=\Un_{i\leq m}\tSG^i[\eseqint{F}{1}{n}],\quad \tSG^{<m}[\eseqint{F}{1}{n}]:=\Un_{i< m}\tSG^i[\eseqint{F}{1}{n}].
\end{gather*}

\bthm
\label{thm:BHvSz}
Assume that $m$ is a natural number and $\roth$ contains a multiplicative idempotent, large for a sequence $\eseq{\cR}\sub\Fin$.
Then for each coloring of  $[\bbN]^m$, there are sets $\eseq{F}\in\Fin$ and elements $y_1\in F_1, y_2\in F_2,\dotsc$ with the following properties.

\be 

\item 
All sequences in the product $F(H_1)\x F(H_2)\x\dotsb$ are proper for all sequences $H_1<H_2<\dotsb$ in $\Fin$.

\item All $m$-sumgraphs of sequences $F(H_1),F(H_2),\dotsc$ have the same color for all sequences $H_1<H_2<\dotsb$ in $\Fin$.
\ee
\ethm

\bpf
Since Theorem~\ref{thm:BHvSz} for $m=2$, implies Theorem~\ref{thm:BHvSz} for $m=1$, assume that $m\geq 2$. 
Let $e$ be a multiplicative idempotent ultrafilter in $\beta\bbN$, large for the sequence $\eseq{\cR}$.
Proceed as in the beginning of the proof of Theorem~\ref{thm:mainfin}, defining a coloring $\chi$ of the set $[\bbN]^{\leq m}$.
Let $G\in e$ be a $\chi$-monochromatic subset of $\bbN$ and assume that the color is green.

\textbf{1st step:} Let 
\[
D_1:=G\cap \Un\cR_1.
\]
The set $D_1$ is in $e$.
Since the ultrafilter $e$ is combinatorially large, there is an element $y_1\in D_1^\star$ such that the set $D_1^\star/y_1$ is in $e$.
Then there is a finite set $F_1\sub D_1^\star$ such that $F_1=R_1\cup \{y_1\}$ for some set $R_1$ in $\cR_1$.
We have $\tSG^{\leq m}[F_1]=F_1\sub D_1$, and thus the set $\tSG^{\leq m}[F_1]$ is green.

\textbf{2nd step:}
Since $F_1\sub D_1^\star$, there is a set $B_1\sub D_1$ in  $e$ such that $F_1+B_1\sub D_1$.
Let
\[
D_2:=B_1\cap D_1^\star/y_1\cap G(\tSG^{\leq m}[F_1])\cap\Un\cR_2.
\]
Since  the set $\tSG^{\leq m}[F_1]$ is green, the set  $G(\tSG^{\leq m}[F_1])$ is in $e$, and thus the set $D_2$ is in $e$, as well.
Since the ultrafilter $e$ is multiplicative, there is an element $y_2\in D_2^\star$ such that the set $D_2^\star/y_2$ is in $e$.
Then there is a finite set $F_2\sub D_2^\star$ such that $F_2=R_2\cup \{y_2\}$ for some set $R_2\in\cR_2$.

The set $\tSG^{\leq m}[F_1,F_2]$ is green:
By the previous step, the set $F_1$ is green, and the set $F_2$ is green as a subset of $D_2$.
Since $F_1+B_1\sub D_1$ and $F_2\sub B_1$, the set $F_1+F_2$ is green.
Since $y_2\in D_2^\star\sub D_1^\star/y_1$, we have $y_1\cdot y_2\in D_1^\star\sub G$.
Thus, the set $\tSG^1[F_1, F_2]=F_1\cup F_2\cup (F_1+F_2)\cup\{y_1\cdot y_2\}$ is green.
Since $F_2\sub G(\tSG^1[F_1])$, the set $\tSG^2[F_1, F_2]$ is green.
Thus, the set $\tSG^{\leq  m}[F_1, F_2]=\tSG^1[F_1, F_2]\cup \tSG^2[F_1, F_2]$ is green.

Fix a natural number $n\geq 2$.
Assume that sets $D_1^\star,\dotsc, D_n^\star$ such that the sets $D_1\spst \dotsb \spst D_n$ are in $e$ and finite sets 
$F_1\sub D_1^\star,\dotsc, F_n\sub D_n^\star$ 
such that $F_i=R_i\cup\{y_i\}$ for some set $R_i\in\cR_i$, the sets $D_i^\star/ y_i$ are in $e$ for all $i\in\{1,\dotsc,n\}$ and the set $\tSG^{\leq m}[\eseqint{F}{1}{n}]$ is green, have been already defined.

\textbf{$\mathbf{(n+1)}$ step:}
Since $F_n\sub D_n^\star$, there is a set $B_n\sub D_n$ in $e$ such that $F_n+B_n\sub D_n$.
Let
\[
D_{n+1}:=B_n\cap D_n^\star/ y_n\cap  G(\tSG^{\leq m}[\eseqint{F}{1}{n}])\cap \Un\cR_n.
\]
Since  the set $\tSG^{< m}[\eseqint{F}{1}{n}]$ is green, the set  $G(\tSG^{< m}[\eseqint{F}{1}{n}])$ is in $e$, and thus the set $D_{n+1}$ is in $e$, as well.
Since the ultrafilter $e$ is multiplicative, there is an element $y_{n+1}\in D_{n+1}^\star$ such that the set $D_{n+1}^\star/ y_{n+1}$ is in $e$.
Then there is a finite set $F_{n+1}\sub D_{n+1}^\star$ such that $F_{n+1}=R_{n+1}\cup\{y_{n+1}\}$ for some set $R_{n+1}\in\cR_{n+1}$.

\bclm\label{clm:tmain}
The set $\tSG^{\leq m}[\eseqint{F}{1}{n+1}]$ is green.
\eclm

\bpf

For a natural number $j$, natural numbers $i_1<\dotsb<i_j\leq n+1$, elements $a_{i_1}\in F_{i_1},\dotsc, a_{i_j}\in F_{i_j}$ and elements $y_{i_1},\dotsc, y_{i_n}$,
the elements $a_{i_1}+\dotsb + a_{i_j}$ and $y_{i_1}\dotsm y_{i_n}$ are in $D_{i_1}$:
By the proof of Claim~\ref{clm:main}, it is enough to show that the assertion is true for $y_{i_1}\dotsm y_{i_n}$.
If $j=1$, then $y_{i_1}\in F_{i_1}\sub D_{i_1}$.
Assume that the claim is true for a natural number $j$ and
take the element $y_{i_1}\dotsm y_{i_{j+1}}$.
By the assumption, we have $y_{i_2}\dotsm y_{i_{j+1}}\in D_{i_2}\sub D_{i_1+1}\sub D_{i_1}^\star/ y_{i_1}$.
Thus, $y_{i_1}\cdot y_{i_2}\dotsm y_{i_{n+1}} \in D_{i_1}^\star\sub  D_{i_1}$.

By the above, the set $\tSG^1[\eseqint{F}{1}{n+1}]$ is green.
Let $A\in \tSG^{\leq m}[\eseqint{F}{1}{n+1}]$ be a set such that $\card{A}>1$.
Then there are a natural number $l$ with $1<l\leq n+1$ and an element $b\in A$ such that $A\sm\{b\}\in\tSG^{<m}[\eseqint{F}{1}{l}]$ and $b\in \tSG^1[\eseqint{F}{l+1}{n+1}]$.
By the above, $b\in D_{l+1}\sub G(\tSG^{<m}[\eseqint{F}{1}{l}])$, and thus the set $A$ is green.
\epf

According to this procedure, define sets $\eseq{D^\star}$ in $e$ and finite sets $F_1\sub D_1^\star, F_2\sub D_2^\star,\dotsc$ with the above properties.

(1)
Fix a sequence $H_1<H_2<\dotsb$ in $\Fin$ and let $\eseq{a}$ be a sequence in $F(H_1)\x F(H_2)\x \dotsb$.
Let $G_1,G_2\in\Fin$ be sets with $G_1< G_2$ and $n=\max G_1$.
Since $a_{G_2}\in \tSG^1[F_{n+1},F_{n+2},\dotsc]$, we have $a_{G_2}\in D_{n+1}$.
Since $a_{G_1}\in\tSG^1[\eseqint{F}{1}{n}]$ and $D_{n+1}\sub G(\tSG^{<m}[\eseqint{F}{1}{n}])$, we have $a_{G_1}\neq a_{G_2}$.

(2) It follows from Claim~\ref{clm:tmain}.
\epf

\section{Comments and open problems}

As we already mentioned in the introduction, Tsaban~\cite{tsramsey} considered theorems about colorings of natural numbers concerning colorings and covers of countable discrete spaces.
This approach can be also applied to theorems considered here.
By Theorem~\ref{thm:mainfin}, we have the following result.

\bprp\label{prp:swb}
Let $X$ be a separable metrizable space and $\cup$ be the semigroup operation on $\tau$.
Assume that $\cB\in\{\Lambda,\Om\}$, the space $X$ satisfies $\sone(\Om,\cB)$, the family $\Om$ contains a free idempotent, large for a sequence $\eseq{\cR}\sub\fin(\tau)$, and $m$ is a natural number.
Then for each coloring of $[\tau]^m$, there are sets $\cF_1\in\cR_1, \cF_2\in\cR_2,\dotsc$ and elements $U_1\in\cF_1, U_2\in\cF_2,\dotsc$ with the following properties.

\be
\item The set $\sset{U_n}{n\in\bbN}$ is in $\cB$.
\item All sequences in the product $\eprod{\cF}$ are proper.
\item The partite $m$-sumgraph of the sequence $\eseq{\cF}$ is monochromatic.
\ee 
\eprp

For a natural number $n$ and a set $F\in\Fin$, we write $n<F$ instead of $\{n\}<F$ and $F<n$ instead of $F<\{n\}$.

\blem\label{lem:corcover}
Let $\eseq{a}$ be a sequence in a semigroup and $e$ be an idempotent ultrafilter on a subsemigroup $S:=\FS(\eseq{a})$, large for a family $\cR\sub\fin(S)$ such that 
\[
\smallmedset{\FS(a_n, a_{n+1},\dotsc)}{n\in\bbN}\sub e\text{ and }\bigcap_{n\in\bbN}
\FS(a_n,a_{n+1},\dotsc)=\emptyset.
\]
Then the family
\[
\cA:= \smallmedset{A\sub\Fin}{
\sset{a_F}{F\in A\text{ and }n<F}\in e\text{ for all }n}
\] 
is an idempotent superfilter on $\Fin$, with the semigroup operation $\cup$, large for the family 
\[
\cR':=\sset{R'\sub\Fin}{\sset{a_F}{F\in R'}\in\cR}.
\]
\elem

\bpf
The family $\cA$ is a superfilter.
Fix a set $A\sub\Fin$ and assume that the set 
\[
B:=\set{G\in\Fin}{\Exists{C\in\cA}\End{\sset{G\cup H}{H\in C}\sub A}}
\]
is in $\cA$.
Let $G\in B$ and $C\in\cA$ be a set such that $\sset{G\cup H}{H\in C}\sub A$.
There is a natural number $n$ such that $G< n$.
Then the set $\sset{H\in C}{n<H}$ is in $\cA$.
Thus, we may assume that $G\cap H=\emptyset$ for all sets $H\in C$.
We have 
\[
a_G + \sset{a_H}{H\in C} = \sset{a_G + a_H}{H\in C}=\sset{a_{G\cup H}}{H\in C}\sub\sset{a_F}{F\in A}
\]
and $\sset{a_H}{H\in C}\in e$.
Since the set $\sset{a_G}{G\in B}$ is in $e$, its superset
\[
\sset{a\in S}{\Exists{D\in e}\End{a + D\sub \sset{a_F}{F\in A}}},
\]
is in $e$, as well.
Since the ultrafilter $e$ is idempotent, we have $\sset{a_F}{F\in A}\in e$.
Fix a natural number $n$.
Since the set $\sset{G\in B}{n< G}$ is in $\cA$, its superset 
\[
\set{G\in\Fin}{\Exists{C\in\cA}\End{\sset{a_{G\cup H}}{H\in C}\sub \sset{a_F}{F\in A\text{ and }n<F}}}
\]
is in $\cA$, as well.
By the above, the set $\sset{a_F}{F\in A\text{ and }n< F}$ is in $e$.
We have $A\in\cA$, and thus the superfilter $\cA$ is idempotent.

Fix a set $A\in\cA$.
Then the set $\sset{a_F}{F\in A}$ is in $e$.
Since the ultrafilter $e$ is large for the family $\cR$, there is a set $R\in \cR$ with $R\sub \sset{a_F}{F\in A}$.
Then there is a finite nonempty set $R'\sub A$ such that $R=\sset{a_F}{F\in R'}$.
We have $R' \in\cR'$, and thus the family $\cA$ is large for the family $\cR'$.
\epf

\bexm
Let $+$ be a semigroup operation on $\bbN$ and $\eseq{a}$ be the sequence $1,2,\dotsc$.
Assume that $\roth$ contains an idempotent, large for a sequence $\eseq{\cR}\sub\Fin$.
Let $m$ be a natural number and fix a coloring $\chi$ of $[\bbN]^m$.
Consider a discrete space $X:=\Fin$ with the semigroup operation $\cup$.
For each natural number $n$, let 
\[
O_n:=\sset{F\in \Fin}{n\notin F}\text{ and }\cR_n':=\sset{A\sub\Fin}{\sset{a_F}{F\in A}\in\cR}.
\]
If $F=\{i_1,\dotsc, i_n\}\in\Fin$, where $i_1<\dotsb<i_n$, then $a_F=i_1+\dotsb +i_n$.
By Lemma~\ref{lem:corcover}, the family 
\[
\cA:= \smallmedset{A\sub\Fin}{
\sset{a_F}{F\in A\text{ and }n<F}\in e\text{ for all }n}
\]
is an idempotent superfilter on $X$, large for the sequence $\eseq{\cR'}$.
By the result of Tsaban~\cite[Theorem~2.7]{tsramsey}, there is a free idempotent in $\cA$.
This idempotent is large for the sequence $\eseq{\cR'}$.
Let $\cup$ be the semigroup operation on $\tau$ and $S:=\FS(\eseq{O})$ be the subsemigroup of $\tau$.

The map from $\Fin$ to $S$ defined by $F\mapsto O_F$ is a semigroup isomorphism.
Fix a set $A\in\cA$.
Then there are sets $\eseq{F}\in A$ with $F_1<F_2<\dotsb$.
By the definition of the sets $O_n$, the family $\sset{O_{F_n}}{n\in\bbN}$ is in $\Om$.
Since 
\[
\sset{O_{F_n}}{n\in\bbN}\sub\sset{O_F}{F\in A},
\]
the later set is in $\Om$, as well.
Thus, 
\[
\cA'':=\smallmedset{\sset{O_F}{F\in A}\sub\Om}{A\in\cA}.
\]
is an idempotent superfilter and it contains a free idempotent, large for the sequence
\[
\sset{\sset{O_F}{F\in\cF}}{\cF\in\cR_1'},
\sset{\sset{O_F}{F\in\cF}}{\cF\in\cR_2'},\dotsc.
\]

Define a coloring $\kappa$ of $[S]^m$ such that for pairwise different sets $H_1, \dotsc ,H_m$ in $\Fin$ we have
\[
\kappa(\{O_{H_1},\dotsc, O_{H_m}\}):=
\begin{cases}
\chi(\{a_{H_1},\dotsc, a_{H_m}\}),&\text{ if }\{a_{H_1},\dotsc, a_{H_m}\}\in[S]^m,\\
1,& \text{ otherwise }.
\end{cases}
\]
By Proposition~\ref{prp:swb}, there are sets $\cF_1\in\cR_1', \cF_2\in\cR_2',\dotsc$ such that all sequences in the product
\[
\sset{O_F}{F\in\cF_1}\x\sset{O_F}{F\in\cF_2}\x\dotsb
\]
are proper and the $m$-sumgraph of the sequence 
$
\sset{O_F}{F\in\cF_1},\sset{O_F}{F\in\cF_2},\dotsc
$
 is $\kappa$-monochromatic.
Let say that the color is green.
We may assume that $\Un\cF_1<\Un\cF_2<\dotsb$.
Let 
\[
R_n:=\sset{a_F}{F\in \cF_n}
\]
for all $n$.

Fix a sequence $\eseq{b}$ in the product $\eprod{R}$.
Then there are sets $F_1\in\cF_1, F_2\in\cF_2,\dotsc$ such that this sequence is equal to the sequence $a_{F_1}, a_{F_2}, \dotsc$.

Let $G=\{i_1,\dotsc,i_n\}$ and $H=\{j_1,\dotsc, j_k\}$ be sets in $\Fin$ with $i_1<\dotsb<i_n<j_1<\dotsb <j_k$ for some natural numbers $n,k$.
We have 
\[
b_G=a_{F_{i_1}}+\dotsb +a_{F_{i_n}}=a_{F_{i_1}\cup\dotsb \cup F_{i_n}}=a_{F_G}
\text{ and }
b_H=a_{F_{j_1}}+\dotsb +a_{F_{j_k}}=a_{F_{j_1}\cup\dotsb \cup F_{j_n}}=a_{F_H}
\]
Since $\{F_{i_1},\dotsc, F_{i_n}\}\in \Un_{i\leq i_n}\cF_i$ and $\Un_{i\leq i_n}\cF_i< \Un_{i=j_1}^{j_k}\cF_i$, we have $b_G\neq b_H$, and thus the sequence $\eseq{b}$ is proper.

Let $H_1<\dotsb<H_m$ be sets in $\Fin$.
Then
\[
b_{H_1}=a_{F_{H_1}},\dotsc, b_{H_m}=a_{F_{H_m}}.
\]
The set $\{a_{F_{H_1}},\dotsc, a_{F_{H_m}}\}$ is in $[S]^m$, and thus
\[
\kappa(\{O_{F_{H_1}},\dotsc, O_{F_{H_m}}\})
=\chi(\{a_{F_{H_1}},\dotsc, a_{F_{H_m}}\}))=
\chi(\{b_{H_1},\dotsc, b_{H_m}\}).
\]
Since $\{O_{F_{H_1}},\dotsc, O_{F_{H_m}}\}$ is an element of the partite $m$-sumgraph of the sequence
$
\sset{O_F}{F\in\cF_1},\sset{O_F}{F\in\cF_2},\dotsc
$, the color of $\{b_{H_1},\dotsc, b_{H_m}\}$ is green.
\eexm

The following problem was formulated by Tsaban~\cite[Subsection~7.3.]{tsramsey}.
Because of its importance we repeat it here.

\bprb[{Tsaban~\cite[Subsection~4.6]{tsramsey}}]
\label{prb:men}
Let $X$ be a space and $\cup$ be the semigroup operation on $\tau$.
Assume that for every family $\cU\in \Lambda$ and coloring of  $[\tau]^2$ there are pairwise disjoint sets $\eseq{\cF}\in\fin(\cU)$ with the following properties.
\be
\item The family $\Un_{n\in\bbN}\cF_n$ is in $\Lambda$.
\item The sequence $\eseq{\Un\cF}$ is proper.
\item The sumgraph of the sequence $\eseq{\Un\cF}$ is monochromatic.
\ee
Is the space $X$ Menger?
\eprb

The property from Problem~\ref{prb:men} is formally weaker than the Menger property.
Thus, the question is whether these two properties are equivalent?

\label{sec:comments}

\section*{Acknowledgments}

I would like to thank Boaz Tsaban, who introduced me to this topic, for encouragement to continue his work and his great impact to my research.
I would like to thank Marion Scheepers, who draw my attention that for a space, if $\cA\cup \cB\in\om{X}$, then one of the families $\cA$ or $\cB$ is in $\om{X}$.


\begin{thebibliography}{99}

\Pa{arhcft}{A. Arhangel'ski\u{\i}}{Hurewicz spaces, analytic sets, and fan tightness of function spaces}{Soviet Mathematics Doklady}{33}{1986}{396}{399}

\Pa{Aurichi}{L. Aurichi}{D-spaces, topological games, and selection principles}{Topology Proceedings}{36}{2010}{107}{122}

\Pa{BaTs}{T. Bartoszy\'nski, B. Tsaban}{Hereditary topological diagonalizations and the Menger--Hurewicz Conjectures}{Proceedings of the American Mathematical Society}{134}{2006}{605}{615}

\Pa{BHcell}{V. Bergelson, N. Hindman}{A combinatorially large cell of a partition on $\bbN$}{Journal of Combinatorial Theory, Series A}{48}{1988}{39}{52}

\Pa{BH}{V. Bergelson, N. Hindman}{Ultrafilters and multidimensional Ramsey theorems}{Combinatorica}{9}{1989}{1}{7}

\Pa{ChRZd}{D. Chodounsk\'{y}, D. Repov\v{s}, L. Zdomskyy}{Mathias forcing and combinatorial covering properties of filters}{The Journal of Symbolic Logic}{80}{2015}{1398}{1410}

\Pa{dde}{W. Deuber}{Partitionen und lineare Gleichungssysteme}{Mathematische Zeitschrift }{133}{1973}{109}{123}

\Pa{DH}{W. Deuber, N. Hindman}{Partitions and sums of $(m,p,c)$-sets}{Journal of Combinatorial Theory, Series A}{45}{1987}{300}{302}

\Pa{E}{R. Ellis}{Distal transformation groups}{Pacific Journal of Mathematics}{8}{1958}{401}{405}

\Pa{FrMill}{D. Fremlin, A. Miller}{On some properties of Hurewicz, Menger and Rothberger}{Fundamenta Mathematicae}{129}{1988}{17}{33}

\Pa{gn}{J. Gerlits, Zs. Nagy}{Some Properties of $\C(X)$, I}{Topology and its Applications}{14}{1982}{151}{161}

\bibitem{Rbook}R. Graham, B. Rothschild, J. Spencer, \textbf{Ramsey Theory}, Wiley, New York ,1980.

\bibitem{Hure25} W.~Hurewicz,
\emph{\"Uber eine Verallgemeinerung des Borelschen Theorems},
Mathematische Zeitschrift \textbf{24} (1925), 401--421.

\bibitem{Hure27} W.~Hurewicz,
\emph{\"Uber Folgen stetiger Funktionen}, Fundamenta Mathematicae \textbf{9} (1927), 193--204.

\Pa{coc2}{W.~Just, A.~Miller, M.~Scheepers, P.~Szeptycki}{The combinatorics of open covers II}{Topology and its Applications}{73}{1996}{241}{266}

\Pa{Hu}{N. Hindman}{Ultratilters and combinatorial number theory}{Lecture Notes in Mathematics}{751}{1979}{119}{184}

\Pa{HFS}{N. Hindman}{Finite sums from sequences within cells of a partition of $\bbN$}{Journal of Combinatorial Theory, Series A}{17}{1974}{1}{11}

\bibitem{HS}N. Hindman, D. Strauss, \emph{Algebra in the Stone--\v{C}ech Compactification}, 2nd ed., de Gruyter, 2012.

\Pa{Laver}{R. Laver}{On the consistency of Borel’s conjecture}{Acta Mathematicae}{137}{1976}{151}{169}

\bibitem{Menger24} K.~Menger,
\emph{Einige \"Uberdeckungss\"atze der Punktmengenlehre},
Sitzungsberichte der Wiener Akademie \textbf{133} (1924), 421--444.


\Pa{mill}{K. Milliken}{Ramsey's theorem with sums or unions}{Journal of Combinatorial Theory, Series A}{18}{1975}{276}{290}

\Pa{N}{K. Numakura}{On bicompact semigroups}{Mathematical Journal of Okayama University}{1}{1952}{99}{108}

\bibitem{Paw94}
J. Pawlikowski,
\emph{Undetermined sets of point-open games},
Fundamenta Mathematicae \textbf{144} (1994), 279--285.

\Pa{rothb}{F. Rothberger}{Eine Versch\"arfung der Eigenschaft C}{Fundamenta Mathematicae}{30}{1938}{50}{55}

\bibitem{sakai} M.~Sakai, \emph{Property C'' and function spaces}, Proceedings of the American Mathematical Society \textbf{104} (1988), 917--919. 

\Pa{schtelg}{M. Scheepers}{A direct proof of a theorem of Telg\'{a}rsky}{Proceedings of the American Mathematical Society}{123}{1995}{3483}{3485} 

\Pa{coc1}{M. Scheepers}{Combinatorics of open covers I: Ramsey theory}{Topology and its  Applications}{69}{1996}{31}{62}

\Pa{coc6}{M. Scheepers}{Combinatorics of open covers, VI. Selectors for sequences of dense sets}{Quaestiones Mathematicae}{22}{1999}{109}{130}

\Pa{schpart}{M. Scheepers}{Open covers and partition relations}{Proceedings of the American Mathematical Society}{127}{1999}{577}{581}

\Pa{Comb}{P. Szewczak, B. Tsaban}{Products of Menger spaces: a combinatorial approach}{Annals of Pure and Applied Logic}{168}{2017}{1}{18}

\bibitem{Mgame}P. Szewczak, B. Tsaban: \emph{Conceptual proofs of the Menger and Rothberger games},
Topology and its Applications \textbf{272} (2020), 107048.

\Pa{pMWien}{P. Szewczak, B. Tsaban, L. Zdomskyy}{Finite powers and products of Menger sets}{
Fundamenta Mathematicae}{253}{2021}{257}{275}

\Pa{tylor}{A. Taylor}{A canonical partition relation for finite subsets of $\w$}{Journal of Combinatorial Theory, Series A}{21}{1976}{137}{146}

\Pa{telg}{R.Telg\'{a}rsky}{On games of Tops\o{}e}{Mathematica Scandinavica}{54}{1984}{170}{176}

\Pa{tsramsey}{B. Tsaban}{Algebra, selections, and additive Ramsey theory}{Fundamenta Mathematicae}{240}{2018}{81}{104}

\Pa{sfh}{B. Tsaban, L. Zdomskyy}{Scales, fields, and a problem of Hurewicz}{Journal of the European Mathematical Society}{10}{2008}{837}{866}

\Pa{vdW}{B. van der Waerden}{Beweiseiner Baudetschen Vermutung}{Nieuw Archief voor Wiskunde
}{19}{1927}{212}{216}

\bibitem{Wiki} Wikipedia, \emph{Selection principle},
\url{https://en.wikipedia.org/wiki/Selection_principle}

\Pa{Mill}{L. Zdomskyy}{Products of Menger spaces in the Miller model}{Advances in Mathematics}{335}{2018}{170}{179}


\end{thebibliography}
\end{document}